\documentclass[11pt]{amsart}
\usepackage{amssymb,mathrsfs,graphicx,subfigure,enumerate}
\usepackage{amsmath,amsfonts,amssymb,amscd,amsthm,bbm}
\usepackage[retainorgcmds]{IEEEtrantools}
\usepackage{colortbl}
\usepackage{tikz}
\usetikzlibrary{patterns}
\usetikzlibrary{positioning}
\usepackage{caption}

\allowdisplaybreaks
\title[Relaxation dynamics of a SIR-flock]{Relaxation dynamics of SIR-flocks with random epidemic states}

\author[Ha]{Seung-Yeal Ha}
\address[Seung-Yeal Ha]{\newline Department of Mathematical Sciences and Research Institute of Mathematics, \newline Seoul National University, Seoul 08826,  Republic of Korea}
\email{syha@snu.ac.kr}

\author[Park]{Hansol Park}
\address[Hansol Park]{\newline Department of Mathematical Sciences\newline Seoul National University, Seoul 08826, Republic of Korea}
\email{hansol960612@snu.ac.kr}

\author[Yang]{Seoyeon Yang}
\address[Seoyeon Yang]{\newline Department of Mathematical Sciences\newline Seoul National University, Seoul 08826, Republic of Korea}
\email{aromeyang@snu.ac.kr}

\newtheorem{theorem}{Theorem}[section]
\newtheorem{lemma}{Lemma}[section]
\newtheorem{corollary}{Corollary}[section]
\newtheorem{proposition}{Proposition}[section]
\newtheorem{remark}{Remark}[section]

\newenvironment{remarks}{{\flushleft {\bf Remarks.}}}{}

\newcommand{\bbr}{\mathbb R}

\makeatletter
\@namedef{subjclassname@2020}{\textup{2020} Mathematics Subject Classification}
\makeatother

\begin{document}

\date{\today}


\thanks{\textbf{Acknowledgment.} The work of S.-Y. Ha is supported by National Research Foundation of Korea (NRF-2017R1A5A1015626).}

\begin{abstract}
We study the collective dynamics of a multi-particle system with three epidemic states as an internal state. For the collective modeling of active particle system, we adopt modeling spirits from the swarmalator model and the SIR epidemic model for the temporal evolution of particles' position and internal states. Under suitable assumptions on system parameters and non-collision property of initial spatial configuration, we show that the proposed model does not admit finite-time collisions so that the standard Cauchy-Lipschitz theory can be applied for the global well-posedness. For the relaxation dynamics, we provide several sufficient frameworks leading to the relaxation dynamics of the proposed model. The proposed sufficient frameworks are formulated in terms of system parameters and initial configuration. Under such sufficient frameworks, we show that the state configuration relaxes to the fixed constant configuration via the exponentially perturbed gradient system and explicit dynamics of the SIR model. We present explicit lower and upper bounds for the minimal and maximal relative distances. 
\end{abstract}

\maketitle \centerline{\date}


\section{Introduction} \label{sec:1}
\setcounter{equation}{0} 
The purpose of this paper is to continue the studies begun in \cite{HJKPZ-0, HJKPZ} on the collective dynamcis modeling of active particles with internal states. Collective behaviors of complex systems are ubiquitous in nature, e.g., crowd dynamics \cite{A-B-G}, aggregation of bacteria \cite{B-B1, B-C-L, B-L-R, K-C-B}, flocking of birds \cite{B-C-C, B-H, C-C-M-P, C-C-G,C-S, H-Liu, H-R, H-T,  M-T1, M-T2, T-T, V-C-B-C-S}, synchronization of fireflies \cite{B-B,Wi} and swarming of fish \cite{D-M}, etc. See survey articles \cite{A-B, A-B-F, B-H-O, D-F-T, H-K-P-Z2, P-R-K, V-Z}. In recent years, thanks to the emerging applications to the decentralized control of multi-particle systems, collective behaviors have received a lot of attention from diverse scientific disciplines such as applied mathematics, biology, control theory and statistical physics etc. In this work, we are interested in the first-order modeling of active particles with random epidemic states (susceptible($S$), infected($I$) and recovered($R$)) as an internal state, i.e., we assume that particle $i$ can take aforementioned three epidemic states with certain probabilities denoted by $S_i, I_i$ and $R_i$ whose precise meaning will be clarified in a minute.  Thus, the state of  particle $i$ is represented by position vector $x_i \in \bbr^d$ and probability vector for epidemic states $W_i = (S_i, I_i, R_i) \in [0, 1]^3$, respectively. In what follows, we briefly describe how to model the dynamics of position vectors and epidemic vectors via continuous dynamical systems. Let $N$ be the size of system, i.e., the total number of particles in a given ensemble $\{ (W_i, x_i) \}_{i=1}^{N}$. 

First, we use the position dynamics $x_i = x_i(t)$ of the swarmalator model \cite{OEK, OHS} which describes the attractive and repulsive forces:
\begin{equation} \label{A-0-0}
 \dot{x}_i =\displaystyle\frac{\kappa_2}{N} \sum_{\substack{1\leq j\leq N\\j\neq i}}  \Psi^{ij}_a \frac{x_j-x_i}{\|x_j-x_i \|^\alpha} + \frac{\kappa_3}{N}\sum_{\substack{1\leq j\leq N\\j\neq i}} \Psi^{ij}_r \frac{x_i -x_j}{\|x_i-x_j\|^\beta}, \quad  i \in [N] := \{1, 2, \cdots, N \}, 
 \end{equation}
 where $\kappa_2, \kappa_3$ are nonnegative coupling strengths, and $\Psi^{ij}_a,~\Psi^{ij}_r$ represent attractive and repulsive weights whose explicit functional forms will be discussed in Section \ref{sec:3}, and we assume that positive system parameters $\alpha, \beta$ satisfy the relation:
 \begin{equation} \label{Assumption-1}
 1\leq \alpha<\beta
 \end{equation}
 so that repulsive force is dominant in a small relative distance regime. Here $\| \cdot \|$ denotes the standard $\ell^2$-norm in ${\mathbb R}^d$. 
 
Next, we use the modeling spirit of the SIR model for the dynamics of epidemic state $W_i$. We introduce convex set ${\mathcal S}$ consisting of all admissible state vectors:
\begin{equation} \label{A-0-1}
\mathcal{S}:= \Big \{W := (S, I, R)\in [0, 1]^3: ~~S + I + R=1 \Big \},
\end{equation}
and we set 
\begin{align*}
\begin{aligned}
& S_i:~\mbox{the probability that the $i$-th particle is in susceptible state}, \\
& I_i:~ \mbox{the probability that the $i$-th particle is in infected state,} \\
& R_i:~\mbox{the probability that the $i$-th particle is in recovered state},
\end{aligned}
\end{align*}
and 
\[ W_i := (S_i, I_i, R_i)~:~\mbox{the epidemic probability vector of the $i$-th particle}. \]
Then, we assume that the dynamics of $W_i$ is governed by the following coupled system:
\begin{equation} \label{A-0-1}
\dot{S}_i =-\kappa_1 \sum_{j=1}^N a^{ij}S_i I_j,\quad \dot{I}_i= \kappa_1 \sum_{j=1}^N a^{ij}S_i I_j-b^i I_i,\quad \dot{R}_i=b^i I_i,~~i \in [N].
\end{equation}
Finally, we couple two systems \eqref{A-0-0} and \eqref{A-0-1} via nonnegative system functions $\Psi^{ij}_a,~\Psi^{ij}_r,~a^{ij}$ and $b^i$ by imposing suitable functional dependences between position variable $\{ x_i \}$ and internal variable $\{ W_i \}$:
\begin{equation} \label{A-0-2}
\Psi_a^{ij} := \Psi_a(W_i, W_j), \quad \Psi^{ij}_r := \Psi_r(W_i, W_j), \quad a^{ij} := a(\|x_i - x_j\|), \quad i, j \in [N].
\end{equation}
Moreover, we also assume that there exist positive constants $\varepsilon_a, \varepsilon_r, M_a$ and $M_r$ such that 
\[  0<\varepsilon_a \leq \Psi_a^{ij} \leq M_a, \quad 
0<\varepsilon_r \leq \Psi_r^{ij} \leq M_r, \quad i, j \in [N]. \]
See Section \ref{sec:3} for explicit functional relations in \eqref{A-0-2}. \newline

\noindent Finally, we combine all three ingredients \eqref{A-0-0}, \eqref{A-0-1} and \eqref{A-0-2} together to write down the dynamical system for $(W_i, x_i)$ with suitable initial data:
\begin{equation}
\begin{cases}  \label{A-1}
\displaystyle \dot{S}_i =-\kappa_1 \sum_{j=1}^N a^{ij}S_i I_j,\quad \dot{I}_i= \kappa_1 \sum_{j=1}^N a^{ij}S_i I_j-b^i I_i,\quad \dot{R}_i=b^i I_i, \quad t > 0,~~i \in [N],  \\
\displaystyle \dot{x}_i =\displaystyle\frac{\kappa_2}{N} \sum_{\substack{1\leq j\leq N\\j\neq i}} \left(\Psi^{ij}_a \cdot\frac{x_j-x_i} {\|x_j-x_i\|^\alpha}\right) + \frac{\kappa_3}{N}\sum_{\substack{1\leq j\leq N\\j\neq i}}\left(\Psi^{ij}_r \cdot\frac{x_i-x_j}{ \|x_i-x_j \|^\beta}\right), \\
\displaystyle  (W_i, x_i)(0)=(W_i^0, x_i^0)  \in {\mathcal S} \times \bbr^d.
\end{cases}
\end{equation}
Throughout the paper, we call the coupled system \eqref{A-1} as the ``{\it SIR-flock model}" for simplicity. 

The goal of this work is to provide sufficient frameworks leading to the emergent dynamics of the SIR-flock model \eqref{A-1}. Now, we briefly discuss four main results. 

First, we show that if spatial configuration is noncollisional initially, there will be no finite-time collisions so that system \eqref{A-1} is globally well-posed by the Cauchy-Lipschitz theory (see Theorem \ref{T2.2}). 

Second, we present a sufficient conditions for the relaxation of epidemic state $W_i$ toward a constant epidemic state. More specifically, our sufficient framework for the relaxation of $W_i$ is expressed in terms of network topology $(a^{ij})$ and a recovering vector $(b^i)$: 
\begin{align}\label{A-2}
a^{ij} :=\begin{cases}
\displaystyle\frac{1}{( \|x_i-x_j\|+L)^\gamma}\quad&\text{if}\quad i\neq j,\\
0&\text{otherwise},
\end{cases}
\hspace{1cm} \mbox{and} \quad   \min_{1 \leq i \leq N} b^i > \frac{\kappa_1(N-1)}{L^\gamma},
\end{align}
where $L$ is a positive constant and $\gamma$ is a nonnegative constant. Under this framework \eqref{A-2},  we show that the epidemic state $W_i$ relaxes to a constant state (see Theorem \ref{T4.1}): there exist a constant state $W_i^{\infty} = (S_i^{\infty}, 0, R_i^{\infty})$ and $\lambda > 0$ such that 
\[
\lim_{t\rightarrow \infty} (S_i(t), I_i(t), R_i(t)) = (S^{\infty}_i, 0, R^{\infty}_i), \quad i \in [N] \quad \mbox{and} \quad  \Big| \frac{1}{N} \sum_{i=1}^N{I_i(t)} \Big|\leq e^{-\lambda t},~~t \geq 0.
\]

Third, we show that if initial configuration and system parameters satisfy suitable conditions, then minimal and maximal relative distances are positive uniformly in time, i.e.,  there exist positive constants $\delta_1$ and $\delta_{\infty}$ such that 
\[
\min_{i \not = j}  \inf_{0\leq t<\infty}  \|x_i(t) - x_j(t) \| \geq\delta_1 \quad \mbox{and} \quad \max_{i \not = j}  \sup_{0\leq t<\infty} \|x_i(t) - x_j(t) \|  \leq \delta_\infty,
\]
see Theorem \ref{T5.1}.

Fourth, we show that the spatial configuration relaxes to a constant configuration asymptotically under a suitable condition on the initial configuration and system parameters (see Theorem \ref{T6.1}). \newline

The rest of this paper is organized as follows. In Section \ref{sec:2}, we briefly review basic properties of the SIR epidemic model and the swarmalator model on the asymptotic relaxation of state variables and present a global well-posedness by verifying nonexistence of finite-time collisions.  In Section \ref{sec:3}, we discuss the modeling spirit of system parameters and coupling functions. In Section \ref{sec:4}, we study the relaxation of epidemic states toward a constant state and in particular, we provide a sufficient condition leading to asymptotic removal of infected particles. In Section \ref{sec:5}, we study the existence of positive lower bound and upper bound for the minimal and maximal relative distances, respectively. In Section \ref{sec:6}, we study a relaxation of spatial configuration toward a fixed spatial configuration using the perturbed gradient flow theory. In Section \ref{sec:7}, we provide several numerical examples and compare them with analytical results obtained in previous sections. Finally, Section \ref{sec:8} is devoted to a brief summary of our main results and some remaining issues for a future work. 

\section{Preliminaries}\label{sec:2}
\setcounter{equation}{0}
In this section, we briefly review basic properties on two related models ``{\it the SIR epidemic model}" and ``{\it the swarmalator model}" which correspond to the subsystems of system \eqref{A-1}. 
\subsection{The SIR epidemic model} \label{sec:2.1}
In this subsection, we briefly discuss the SIR model \cite{KKJ, Va} which is a prototype model for the spread of disease or virus. Kermack and McKendrick's work in 1927 motivated a large number of modelings in epidemics \cite{KM}. In particular, there are some of the works for diffusive SIR models on stability and numerical methods \cite{CDT, SW, ZZHZ}. First, we introduce the following three observables:
\begin{align*}
S& = S(t): \text{ratio of {\it susceptible} individuals in a total population},\\
I &= I(t): \text{ratio of  {\it infected} individuals  in a total population},\\
R &= R(t): \text{ratio of {\it recovered} individuals in a total population}.
\end{align*} 
Then, we assume that the state $(S, I, R)$ is governed by the Cauchy problem to the coupled system of ODEs:
\begin{align}\label{B-1}
\begin{cases}
\displaystyle \dot{S} =-a S I, \quad t > 0, \\
\displaystyle \dot{I} =a S I-b I ,\\
\displaystyle \dot{R} =b I, \\
\displaystyle  (S, I, R)(0) = (S^0, I^0, R^0),
\end{cases}
\end{align}
where $a$ and $b$ are positive constants. This model was designed to explain the spreading of epidemic diseses. We can express the SIR model via the following pictorial diagram. 
\begin{center}
\begin{tikzpicture}
  \node (S) {$S$};
  \node[right=2cm] (I) {$I$};
  \node[right=4cm] (R) {$R$};
  \draw[-stealth,shorten >= 2pt] (S.north) to[bend left] node[midway,above] {$aSI$} (I.north);
  \draw[-stealth] (I.north) to[bend left] node[midway,above] {$bI$} (R.north);
\end{tikzpicture}
\end{center}
Next, we list basic properties of the SIR model in the following lemma.
\begin{proposition} \label{P2.1}
Let $(S, I, R)$ be a global solution of \eqref{B-1} with initial data satisfying
\[  S^0 \geq 0, \quad I^0 \geq 0, \quad R^0 \geq 0 \quad \mbox{and} \quad S^0+I^0+R^0=1. \]
Then, one has 
\[
S(t)+I(t)+R(t)=1,\quad S(t) \geq 0,~~I(t) \geq 0,~~R(t)\geq0, \quad t \geq 0 \quad \mbox{and} \quad \lim_{t\to +\infty}I(t)=0.
\]
\end{proposition}
\begin{proof} (i)~We add three equations in \eqref{B-1} to get 
\[
\dot{S}+\dot{I}+\dot{R}=0.
\]
This yields
\begin{equation} \label{B-2}
S(t) +I(t)+R(t)= S^0 +I^0 +R^0=1 \quad \mbox{for all $t \geq 0$}. 
\end{equation}

\vspace{0.2cm}

\noindent  (ii)~We integrate the first two equations in \eqref{B-1} to obtain
\begin{align}
\begin{aligned} \label{B-3}
S(t) &=S^0 \exp \Big[  -a \int_0^t I(\tau) d\tau \Big ] \geq 0, \\
I(t) &=I^0 \exp \Big[ \int_0^t \Big(aS(\tau)-b \Big)d\tau \Big] \geq 0.
\end{aligned}
\end{align}
The nonnegativity of $R$ follows from the third equation and the nonnegativity of $\eqref{B-3}_2$:
\[ {\dot R} = b I \geq 0. \]
On the other hand, it follows from \eqref{B-2} and \eqref{B-3} that $R$ is bounded by $1$:
\[ R(t)=1-S(t)-I(t)\leq 1. \]
Since $\dot{R}\geq 0$ and is bounded above by 1, there exists an asymptotic constant state $R^{\infty}$ such that 
\begin{equation} \label{B-4}
\exists~R^{\infty} :=\lim_{t\rightarrow\infty}R(t)   \quad \mbox{and} \quad \lim_{t\rightarrow\infty}\dot{R}(t)=0.
\end{equation}

\vspace{0.2cm}

\noindent (iii)~By $\eqref{B-4}$, one has 
\[
\lim_{t \to \infty} I(t) = \frac{1}{b} \lim_{t \to \infty} \dot{R}(t)=0.
\]
\end{proof}
\subsection{The swarmalator model}  Let $x_i$ and $\theta_i$ be the position and phase of the $i$-th swarmalator, respectively. Then, its dynamics is governed by the Cauchy problem to the Kuramoto-type swarmalator model:
\begin{equation} 
\begin{cases}\label{B-5}
\displaystyle {\dot x}_i =w_i+\frac{1}{N}  \sum_{\substack{1\leq j\leq N\\j\neq i}}  \left[ \tilde{\Psi}_a(\theta_j-\theta_i)\frac{x_j-x_i}{\|x_j-x_i\|^\alpha}-\tilde{\Psi}_r(\theta_j-\theta_i)\frac{x_j-x_i}{\|x_j-x_i\|^\beta}\right],\quad t>0,\\
\displaystyle {\dot \theta}_i =\nu_i+\frac{\kappa}{N}  \sum_{\substack{1\leq j\leq N\\j\neq i}}  \frac{\sin(\theta_j-\theta_i)}{\|x_j-x_i\|^\gamma},\quad i \in [N],\\
\displaystyle (x_i(0),\theta_i(0))=(x_{i}^0,\theta_{i}^0).
\end{cases}
\end{equation}
Here, constants $w_i$ and $\nu_i$ are the natural velocity and frequency of the $i$-th particle, respectively,  $\| \cdot \|$ denotes the standard Euclidean $\ell^2$-norm in $\mathbb{R}^d$. System parameters and interaction functions  ${\tilde \Psi}_a$ and ${\tilde \Psi}_r$ satisfy 
\begin{align}
\begin{aligned}\label{B-6}
& 1\leq \alpha<\beta, \quad {\tilde \Psi}_a(\theta)= {\tilde \Psi}_a(-\theta),\quad \tilde{\Psi}_r(\theta)= \tilde{\Psi}_r(-\theta),\quad \theta\in\mathbb{R},\\
& 0< \tilde{m}_a\leq \tilde{\Psi}_a(\theta)\leq \tilde{M}_a<\infty,\quad 0< \tilde{m}_r\leq \tilde{\Psi}_r(\theta)\leq \tilde{M}_r<\infty.
\end{aligned}
\end{align}
If we simply set $\gamma=0$, system $\eqref{B-5}_2$ becomes the Kuramoto model \cite{Ku2}:
\begin{equation*} \label{B-7}
{\dot \theta}_i=\nu_i+\frac{\kappa}{N}\sum_{j =1}^{N}{\sin(\theta_j-\theta_i)}.
\end{equation*}
Due to the singular terms in the swarmalator model $\eqref{B-5}_1$, it is important to make sure that there is no finite-time collision which can be quantified in the following proposition.  \newline

First, we set several parameters:
\begin{align}
\begin{aligned} \label{B-7}
& \gamma_1:=\frac{2\tilde{m}_a}{N},\quad \gamma_2:=\frac{2\tilde{M}_r}{N},\quad \gamma_3:=D(W) := \max_{i,j} |w_i - w_j|, \\
& {\mathcal Q} := \frac{N (N-1)}{2}, \quad  \tilde{B}:= \frac{\tilde{m}_r}{2(\tilde{M}_a + \tilde{M}_r) + {\mathcal D}(W)},  \\
& D_1(t):= \min_{1 \leq i \not =  j \leq N} \| x_i(t) - x_j(t) \|, \quad   D(t):= \max_{1 \leq i \not =  j \leq N} \| x_i(t) - x_j(t) \|,   \\
& \Lambda := \Big( \frac{1}{2} \Big)^{{\mathcal Q}} \max \{  D_1(0), 1 \}^{\mathcal Q} \Big(  \Big( \max\Big\{ \Big( \frac{\tilde{m}_r}{ 2\tilde{M}_a} \Big)^{\frac{1}{\beta -\alpha}},   1\Big \} \Big)^{\mathcal Q} \Big( \frac{\tilde{B}}{N} \Big )^{\frac{Q}{\beta - 1}}.
\end{aligned}
\end{align}
Next, we state two results on the positivity of minimal  and maximal distance between particles.

\begin{proposition} \label{P2.2}
\emph{(Positivity of minimal distance) \cite{HJKPZ}}
Suppose system parameters satisfy \eqref{B-6}, and the initial data $(X^0,\Theta^0)$ is non-collisional:
\[
\kappa > 0, \qquad \min_{1 \leq i \not = j \leq N} \|x^{0}_{i}-x^{0}_{j} \|>0.
\]
Then, there exists a global solution $(X,\Theta)$ to \eqref{B-5} - \eqref{B-6} and $\tilde{\delta}_1 > 0$ such that 
\[
\inf_{0\leq t<\infty}\min_{1 \leq i \not = j \leq N } \|x_i(t)-x_j(t) \|\geq \tilde{\delta}_1.
\]
\end{proposition}

\begin{proposition} \label{P2.3}
\emph{\cite{HJKPZ}}
Suppose system parameters and initial data satisfy the following framework:
\begin{itemize}
\item
~$({\mathcal F}1)$ : For $\alpha=1$,
\begin{equation} \label{B-7-1}
\Lambda >\left(\frac{\tilde{M}_r}{\tilde{m}_a}\right)^{\beta-\alpha},\quad \gamma_1>\gamma_3.
\end{equation}
\item
~$({\mathcal F}2)$ : For $\alpha>1$,
\begin{equation} \label{B-7-2}
\begin{cases}
\displaystyle D(X^0)<y^*,\quad \Lambda >\left(\frac{\tilde{M}_r}{\tilde{m}_a}\right)^{\beta-\alpha}, \\
\displaystyle -\gamma_1\left(\frac{\gamma_1(\alpha-1)}{\gamma_2(\beta-1)}\right)^{(\beta-\alpha)(\alpha-1)}+\gamma_2\left(\frac{\gamma_1(\alpha-1)}{\gamma_2(\beta-1)}\right)^{(\beta-\alpha)(\beta-1)}+\gamma_3\leq 0,
\end{cases}
\end{equation}
where $y^*$ is the largest root of the equation:
\[
-\gamma_1 y^{1-\alpha}+\gamma_2 y^{1-\beta}+\gamma_3=0,
\]
\end{itemize}
and let $(X,\Theta)$ be a solution to \eqref{B-5} - \eqref{B-6}. Then, there exists a positive constant $\tilde{\delta}_{\infty}$ such that 
\[
\sup_{0\leq t<\infty}D(t)\leq \tilde{\delta}_{\infty}.
\]
\end{proposition}
\begin{proof}
For a proof, we refer to \cite{HJKPZ}.
\end{proof}

\vspace{0.5cm}

Consider the following ansatz for $\tilde{\Psi}_a$ and $\tilde{\Psi}_r$:
\[
\tilde{\Psi}_a(\theta) :=1+J\cos\theta,\quad \tilde{\Psi}_r :=1,\quad |J|<1.
\]
In this case, system \eqref{B-5} becomes 
\begin{equation} \label{B-8}
\begin{cases}
 \displaystyle   {\dot x}_i =w_i+\frac{1}{N}  \sum_{\substack{1\leq j\leq N\\j\neq i}}  \left[(1+J\cos(\theta_j-\theta_i))\frac{x_j-x_i}{\|x_j-x_i\|^\alpha}-\frac{x_j-x_i}{\|x_j-x_i\|^\beta}\right],\quad t>0,\\
 \displaystyle  {\dot \theta}_i  =\nu_i+\frac{\kappa}{N}  \sum_{\substack{1\leq j\leq N\\j\neq i}} \frac{\sin(\theta_j-\theta_i)}{\|x_j-x_i\|^\gamma}.
\end{cases}
\end{equation}
By  Proposition \ref{P2.1} and  convergence result of the perturbed gradient system, system \eqref{B-8} exhibits phase synchronization. 
\begin{theorem} \label{T2.1}
\emph{\cite{HJKPZ}}
Suppose system parameters, natural velocities and initial data satisfy one of the frameworks \eqref{B-7-1} and \eqref{B-7-2}, and 
\[ \kappa>0, \quad D(\nu) =0 \quad \mbox{and} \quad  0<D(\Theta^0)<\pi, \]
and let $(X,\Theta)$ be a solution to \eqref{B-8}. Then there exists asymptotic state $X^\infty$ such that
\[
\lim_{t\rightarrow \infty}X(t)=X^\infty.
\]
\end{theorem}
\begin{proof}
For a proof, we refer to Section 4.2 of \cite{HJKPZ}.
\end{proof}



\subsection{A global well-posedness} \label{sec:2.3}
In this subsection, we discuss a global well-posedness of system \eqref{A-1}. Since the R.H.S. of $\eqref{A-1}_2$ contains the term  $\|x_j - x_i \|$ in the denominators, as long as we can rule out the possibility of finite-time collisions, we obtain a global well-posedness using the standard  Cauchy-Lipschitz theory. In the sequel, we show that system \eqref{A-1} does not admit a finite-time collision, as long as there is no collisions initially. This will be done using a contradiction argument and Gronwall's inequality. \newline

Suppose that initial spatial configuration satisfy
\[ \min_{1 \leq i \not = j \leq N}  \|x_i^0 -x_j^0 \| > 0. \]
Then, by the continuity of solution, there will be no collisions between particles at least small time interval $[0, \varepsilon)$ with $\varepsilon \ll 1$. Then, using the Cauchy-Lipschitz theory, we can show that system \eqref{A-1} has a local smooth solution $\{ (W_i, x_i) \}$ in the time-interval $[0, \varepsilon)$. Now, in order to show a global well-posedness, it suffices to show that there will be no finite-time collisions. Suppose there is a finite-time collision and let $t_0$ be {\it the first collision time}. We take one of the particles that make a collision and fix it by $p$. Then, we define a set containing all the particles involved in the collision at time $t_0$ by ${\mathcal C}$. Now, we set the following handy notation:
\[
{\mathcal R} :=\{1,2,\cdots,N\} \setminus {\mathcal C}, \qquad \chi_{\mathcal C} := \sqrt{\sum_{i,j\in {\mathcal C}} \|x_i-x_j\|^2}, \qquad  \sum_{i \not = j} :=  \sum_{\substack{1\leq j\leq N\\j\neq i}}. 
\]
Then, it is easy to see that 
\begin{equation} \label{B-8-0-0}
 \chi_{\mathcal C}(t_0) = 0, \quad \chi_{\mathcal C}(t) > 0, \quad t \in (0, t_0).
 \end{equation}
By direct calculations, one has
\begin{align}
\begin{aligned} \label{B-8-0}
\frac{d\chi_{\mathcal C}^2}{dt} &= 2\sum_{i,j \in {\mathcal C}}{(x_i -x_j)\cdot (\dot{x_i}-\dot{x_j})}\\
&= 2\sum_{i,j \in {\mathcal C}}{(x_i -x_j)\cdot \left[\frac{\kappa_2}{N}\sum_{k\neq i}{\left(\Psi_a^{ik} \cdot\frac{x_k-x_i}{\|x_k-x_i\|^\alpha}\right)-\frac{\kappa_3}{N}\sum_{k\neq i}\left(\Psi_r^{ik} \cdot\frac{x_k-x_i}{\|x_k-x_i \|^\beta}\right)}\right]} \\
&-2\sum_{i,j \in {\mathcal C}}{(x_i -x_j)\cdot \left[\frac{\kappa_2}{N}\sum_{k\neq j}{\left(\Psi_a^{jk} \cdot\frac{x_k-x_j}{\|x_k-x_j\|^\alpha}\right)-\frac{\kappa_3}{N}\sum_{k\neq j}\left(\Psi_r^{jk} \cdot\frac{x_k-x_j}{\|x_k-x_j\|^\beta}\right)}\right]}\\
&= 2\sum_{i,j\in {\mathcal C}}{(x_i -x_j)\cdot \left[ \frac{\kappa_2}{N}\sum_{k\in C, k\neq i}{\left(\Psi_a^{ik} \cdot\frac{x_k-x_i}{\|x_k-x_i\|^\alpha}\right)}
- \frac{\kappa_2}{N}\sum_{k\in C, k\neq j}{\left(\Psi_a^{ik} \cdot\frac{x_k-x_j}{\|x_k-x_j\|^\alpha}\right)}\right]}\\
&+ 2\sum_{i,j\in {\mathcal C}}{(x_i -x_j)\cdot \left[\frac{\kappa_2}{N}\sum_{k\in R}{\left(\Psi_a^{ik} \cdot\frac{x_k-x_i}{|x_k-x_i|^\alpha}\right)}
- \frac{\kappa_2}{N}\sum_{k\in R}{\left(\Psi_a^{ik} \cdot\frac{x_k-x_j}{\|x_k-x_j\|^\alpha}\right)}\right]}\\
&- 2\sum_{i,j\in {\mathcal C}}{(x_i -x_j)\cdot \left[ \frac{\kappa_3}{N}\sum_{k\in C, k\neq i}{\left(\Psi_r^{ik}\cdot\frac{x_k-x_i}{\|x_k-x_i\|^\beta}\right)}
- \frac{\kappa_3}{N}\sum_{k\in C, k\neq j}{\left(\Psi_r^{ik}\cdot\frac{x_k-x_j}{\|x_k-x_j\|^\beta}\right)}\right]}\\
&- 2\sum_{i,j\in {\mathcal C}}{(x_i -x_j)\cdot \left[\frac{\kappa_3}{N}\sum_{k\in R}{\left(\Psi_r^{ik} \cdot\frac{x_k-x_i}{\|x_k-x_i\|^\beta}\right)}
- \frac{\kappa_3}{N}\sum_{k\in R}{\left(\Psi_r^{ik}\cdot\frac{x_k-x_j}{\|x_k-x_j\|^\beta}\right)}\right]} \\
&=: \mathcal{I}_{11}+\mathcal{I}_{12}+\mathcal{I}_{13}+\mathcal{I}_{14}.
\end{aligned}
\end{align}

\vspace{0.5cm}

In the following lemma, we provide estimates for ${\mathcal I}_{1i}$. 
\begin{lemma} \label{L2.1}
The term ${\mathcal I}_{1i}$ with $1 \leq i \leq 4$ satisfies
\begin{align*}
\begin{aligned}
& (i)~\mathcal{I}_{11} +\mathcal{I}_{13} \geq -\frac{2M_a| {\mathcal C}|}{N}\sum_{\substack{i,k\in {\mathcal C} \\k\neq i}}{\|x_k-x_i\|^{2-\alpha}}+\frac{2m_r|{\mathcal C}|}{N}\sum_{\substack{i,k\in {\mathcal C}\\k\neq i}}{\|x_k-x_i\|^{2-\beta}}, \\
& (ii)~{\mathcal I}_{12} + {\mathcal I}_{14} \geq -\frac{4| {\mathcal R} | \cdot |{\mathcal C}|}{N}\left(\frac{\kappa_2 M_a}{\delta^{\alpha-1}} +\frac{\kappa_3 M_r}{\delta^{\beta-1}} \right)\chi_{\mathcal C},
\end{aligned}
\end{align*}
where $|A|$ is the cardinality of the set $A$ and $\displaystyle \delta := \inf_{\substack{i \in {\mathcal C}, \\ j \in {\mathcal R}}} \| x_i - x_j \|$.
\end{lemma}
\begin{proof}
We provide estimate for $\mathcal{I}_{1i}$ ($i=1,2,3,4$) as follows:
\vspace{0.5cm}

\noindent (i)~(Estimate of $\mathcal{I}_{11}+\mathcal{I}_{13}$):
Note that
\begin{align*}
& 2\sum_{i,j \in {\mathcal C}}{(x_i -x_j)\cdot \left[\frac{\kappa_2}{N}\sum_{k\in {\mathcal C}, k\neq i}{\left(\Psi^{ik}_a \cdot\frac{x_k-x_i}{\|x_k-x_i\|^\alpha}\right)} \right]} \\
& \hspace{1cm}  = \frac{2\kappa_2}{N}\sum_{i,j \in {\mathcal C}}\sum_{\substack{k\in {\mathcal C}\\k\neq i}}{\left(\Psi_a^{ik} \cdot\frac{(x_i -x_j)\cdot (x_k-x_i)}{\|x_k-x_i\|^\alpha}\right)} = \frac{2\kappa_2}{N}\sum_{\substack{i,j,k\in {\mathcal C} \\k\neq i}}{\left(\Psi_a^{ik} \cdot\frac{(x_i -x_j)\cdot (x_k-x_i)}{\|x_k-x_i\|^\alpha}\right)}\\
& \hspace{1cm}  = \frac{\kappa_2}{N}\sum_{\substack{i,j,k\in {\mathcal C} \\k\neq i}}{\left(\Psi_a^{ik} \cdot\frac{(x_i -x_j)\cdot (x_k-x_i)}{\|x_k-x_i\|^\alpha} + \Psi_a^{ik} \cdot\frac{(x_k -x_j)\cdot (x_i-x_k)}{\|x_i-x_k\|^\alpha}\right)}\\
& \hspace{1cm} = -\frac{\kappa_2}{N}\sum_{\substack{i,j,k\in {\mathcal C}\\k\neq i}}{\left( \frac{\Psi_a^{ik}}{\|x_k-x_i\|^{\alpha-2}}\right)} = -\frac{|{\mathcal C}|\kappa_2}{N}\sum_{\substack{i,k\in {\mathcal C}\\k\neq i}}{\left( \frac{\Psi_a^{ik}}{\|x_k-x_i\|^{\alpha-2}}\right)}.
\end{align*}
Therefore, one has
\begin{equation} \label{B-9}
\mathcal{I}_{11} = -\frac{2|{\mathcal C}|\kappa_2}{N}\sum_{\substack{i,k\in {\mathcal C} \\k\neq i}}{\left( \frac{\Psi_a^{ik}}{\|x_k-x_i\|^{\alpha-2}}\right)}.
\end{equation}
Similarly, one has 
\begin{equation} \label{B-10}
\mathcal{I}_{13} = \frac{2|{\mathcal C}|\kappa_3}{N}\sum_{\substack{i,k\in {\mathcal C}\\k\neq i}}{\left( \frac{\Psi_r^{ik}}{\|x_k-x_i\|^{\beta-2}}\right)}.
\end{equation}
We combine \eqref{B-9} and \eqref{B-10} to obtain
\begin{align*}
 \mathcal{I}_{11} +\mathcal{I}_{13}  
&= -\frac{2|{\mathcal C}|\kappa_2}{N}\sum_{\substack{i,k\in {\mathcal C} \\k\neq i}}{\left( \frac{\Psi_a^{ik} }{\|x_k-x_i\|^{\alpha-2}}\right)} + \frac{2|{\mathcal C}|\kappa_3}{N}\sum_{\substack{i,k\in {\mathcal C} \\k\neq i}}{\left(\frac{\Psi_r^{ik}}{\|x_k-x_i\|^{\beta-2}}\right)}\\
&\geq -\frac{2M_a|{\mathcal C}|}{N}\sum_{\substack{i,k\in {\mathcal C} \\k\neq i}}{\|x_k-x_i\|^{2-\alpha}}+\frac{2m_r|{\mathcal C}|}{N}\sum_{\substack{i,k\in {\mathcal C}\\k\neq i}}{\|x_k-x_i\|^{2-\beta}}.
\end{align*}

\vspace{0.5cm}

\noindent (ii)~(Estimate of $\mathcal{I}_{12}+\mathcal{I}_{14}$): Since $|x_k-x_i|\geq\delta$, one has
\begin{align*}
{\mathcal I}_{12} &= 2\sum_{i,j\in {\mathcal C}}{(x_i -x_j)\cdot \left[\frac{\kappa_2}{N}\sum_{k\in R}{\left(\Psi_a^{ik} \cdot\frac{x_k-x_i}{\|x_k-x_i\|^\alpha}\right)}
- \frac{\kappa_2}{N}\sum_{k\in {\mathcal R}}{\left(\Psi_a^{ik} \cdot\frac{x_k-x_j}{\|x_k-x_j\|^\alpha}\right)}\right]}\\
&\geq -\frac{4\kappa_2}{N}\sum_{\substack{i,j\in {\mathcal C} \\k\in {\mathcal R}}}{M_a\frac{\|x_i -x_j\|}{\delta^{\alpha-1}}} = -\frac{4\kappa_2 |{\mathcal R}|M_a}{N \delta^{\alpha-1}}\sum_{i,j\in {\mathcal C}} \|x_i -x_j \|\geq -\frac{4\kappa_2|{\mathcal R}| \cdot |{\mathcal C}|M_a}{N\delta^{\alpha-1}} \chi_{\mathcal C},
\end{align*}
where we used the Cauchy-Schwarz inequality to see
\[
\sum_{i,j\in {\mathcal C}}\|x_i -x_j\| \leq |{\mathcal C}|\chi_{\mathcal C}.
\]
Similarly, one has 
\begin{align*}
{\mathcal I}_{12}  + {\mathcal I}_{14}  \geq -\frac{4|{\mathcal R}| \cdot |{\mathcal C}|}{N}\left(\frac{\kappa_2 M_a}{\delta^{\alpha-1}} +\frac{\kappa_3 M_r}{\delta^{\beta-1}} \right)\chi_{\mathcal C}.
\end{align*}
\end{proof}
It follows from \eqref{B-8-0} and Lemma \ref{L2.1} that 
\begin{align}
\begin{aligned} \label{B-11}
\frac{d\chi_C^2}{dt} &\geq  -\frac{2M_a| {\mathcal C}|}{N}\sum_{\substack{i,k\in {\mathcal C} \\k\neq i}}{\|x_k-x_i\|^{2-\alpha}}+\frac{2m_r|{\mathcal C}|}{N}\sum_{\substack{i,k\in {\mathcal C}\\k\neq i}}{\|x_k-x_i\|^{2-\beta}} \\
&\hspace{0.5cm}  -\frac{4| {\mathcal R} | \cdot |{\mathcal C}|}{N}\left(\frac{\kappa_2 M_a}{\delta^{\alpha-1}} +\frac{\kappa_3 M_r}{\delta^{\beta-1}} \right)\chi_{\mathcal C} \\
&= -c_1\sum_{\substack{i,k\in C\\k\neq i}}{\|x_k-x_i\|^{2-\alpha}}+c_2\sum_{\substack{i,k\in C\\k\neq i}}{\|x_k-x_i\|^{2-\beta}} -c_3\chi_C,
\end{aligned}
\end{align}
where $c_1, c_2, \text{and } c_3$ are positive constants defined as follows.
\[
c_1 := \frac{2M_a| {\mathcal C}|}{N}, \quad c_2 := \frac{2m_r|{\mathcal C}|}{N}, \quad c_3 := \frac{4| {\mathcal R} | \cdot |{\mathcal C}|}{N}\left(\frac{\kappa_2 M_a}{\delta^{\alpha-1}} +\frac{\kappa_3 M_r}{\delta^{\beta-1}} \right).
\]
On the other hand, since $1 \leq \alpha < \beta$ in \eqref{Assumption-1}, there exists a small positive constant $\varepsilon \ll 1$ such that for $t\in(t_0-\varepsilon,t_0)$,
\begin{equation} \label{B-12}
c_1\sum_{\substack{i,k\in C\\k\neq i}}{\|x_k-x_i\|^{2-\alpha}} < \frac{c_2}{4}\sum_{\substack{i,k\in C\\k\neq i}}{\|x_k-x_i\|^{2-\beta}},\quad c_3\chi_C<\frac{c_2}{4}\sum_{\substack{i,k\in C\\k\neq i}}{\|x_k-x_i \|^{2-\beta}}.
\end{equation}
Therefore, we combine \eqref{B-11} and \eqref{B-12} to find 
\[
\frac{d\chi_C^2}{dt} >  \frac{c_2}{2}\sum_{\substack{i,k\in C\\k\neq i}}{ \|x_k-x_i \|^{2-\beta}}, \quad t\in(t_0-\varepsilon,t_0).
\]
Moreover, there exists a positive constant $\varepsilon^* \in (0, \varepsilon)$ such that for $t\in(t_0-\varepsilon^*,t_0)$,
\[
\|x_i -x_j \|^{2-\beta} \geq \|x_i -x_j \|^2, \: \text{for }t\in(t_0-\varepsilon^*,t_0).
\]
Thus, we have
\[
\frac{d}{dt}\chi_C^2 > \frac{c_2}{2}\sum_{\substack{i,k\in C\\k\neq i}}{\|x_k-x_i\|^{2-\beta}} \geq \frac{c_2}{2} \chi_C^2, \quad \text{  for }t\in(t_0-\varepsilon^*,t_0).
\]
i.e.,
\[ \frac{d}{dt}\chi_C^2 >  \frac{c_2}{2} \chi_C^2, \quad \text{  for }t\in(t_0-\varepsilon^*,t_0). \]
We integrate the above differential inequality from $t_0-\varepsilon^*$ to  $t_0$ to get 
\[
e^{-\frac{c_2}{2}t_0}\chi_C^2(t_0) \geq e^{-\frac{c_2}{2}(t_0-\varepsilon^*)}\chi_C^2(t_0-\varepsilon^*).
\]
On the other hand, it follows from \eqref{B-8-0-0} that the left-hand side is zero, but the right-hand side is strictly positive, which is contradictory. Therefore, there are no finite-time collisions between particles from the well-prepared initial configuration. Finally, we can summarize the previous argument as follows. 
\begin{theorem}\label{T2.2}
\emph{(A global well-posedness)}
Suppose that initial spatial configuration is noncollisional in the sense that 
\[ \min_{1 \leq i \not = j \leq N}  \|x_i^0 -x_j^0 \| > 0. \]
Then, there exists a unique global solution $\{ (W_i,x_i) \}$ to system \eqref{A-1} in any finite-time interval. 
\end{theorem}

\section{Modeling of system functions}\label{sec:3}
\setcounter{equation}{0}
In this section, we discuss meaning of the system parameters and system functions appearing in the SIR-flock \eqref{A-1}:
\[ (a^{ij}):~\mbox{interaction topology}, \quad (b^i):~\mbox{recovering vector}, \quad  \Psi_a,~~\Psi_r:~\mbox{coupling  weight functions}.  \]
\subsection{Interaction topology and recovering vector} \label{sec:3.1} For the modeling purpose, we assume that the nonnegative value $a^{ij} = a^{ij}(\|x_i(t)-x_j(t) \|)$ depends on the relative distance $\| x_i - x_j \|$ between the $i$-th particle and the $j$-th particle. To be definiteness, we set
\begin{align}\label{C-1}
a^{ij}=\begin{cases}
\displaystyle\frac{1}{ (\|x_i-x_j\|+L)^\gamma}\quad&\text{if}\quad i\neq j,\\
0&\text{otherwise},
\end{cases}
\end{align}
where $L$ is a positive constant and $\gamma$ is a nonnegative constant. Since the $i$-th particle can not be affected by itself, we assume the diagonal entries of $(a^{ij})$ are zero: 
\[ a^{ii}=0, \quad \forall~i \in [N]. \]

Now we discuss the natural recovering vector $b = (b^1, \cdots, b^N)$. Suppose that every particles have own immune system, the disease will disappear automatically. Thus, we assume
\[ b^i>0, \quad \forall~i \in [N]. \]
If the immune system of the $i$-th particle is well-functioning, then $b^i$ has a large value. In contrast, if the immune system of the $i$-th particle is not well-functioning, then $b^impimwwwi$ will have a small value. In this work, we assume that the natural recovering vector $b_i$ is a positive constant (see the following figure in the sequel).
\newline

\begin{center}
\begin{figure}[h] \label{Fig-A}
\begin{tikzpicture}
\node (Si) at (0, 2) {$S_i$};
\node (Ii) at (3, 2) {$I_i$};
\node (Ri) at (6, 2) {$R_i$};
\node (Sj) at (0, 0) {$S_j$};
\node (Ij) at (3, 0) {$I_j$};
\node (Rj) at (6, 0) {$R_j$};
\draw [->] (Si) -- (Ii);
\draw [->] (Ii) -- (Ri);
\draw [->] (Sj) -- (Ij);
\draw [->] (Ij) -- (Rj);
\node  (A) at (1.5, 2-0.3) {$a^{ij}S_iI_j$};
\node  (B) at (1.5+3, 2-0.3) {$b^iI_i$};
\node (C) at (1.5, 0.3) {$a^{ji}S_jI_i$};
\node (D) at (1.5+3, 0.3) {$b^jI_j$};
\draw [dashed, <-] (A) -- (Ij);
\draw [dashed, <-] (C) -- (Ii);
\end{tikzpicture}
\end{figure}

\end{center}
\subsection{Coupling weight functions} \label{sec:3.2}
Recall that the dynamics of $x_i$ is governed by the following system:
\begin{equation} \label{C-2}
\dot{x}_i = \frac{\kappa_2}{N}  \sum_{j\neq i} \left(\Psi_a^{ij} \frac{x_j-x_i}{ \|x_j-x_i\|^\alpha}\right)-\frac{\kappa_3}{N}\sum_{j\neq i}\left(\Psi_r^{ij} \frac{x_j-x_i}{ \|x_j-x_i\|^\beta}\right).
\end{equation}
Now, the matter of question is how to model $\Psi_a$ and $\Psi_r$ in terms of $W_i'$s.  \newline

If two state vectors $W_i$ and $W_j$ are similar, then particles $x_i$ and $x_j$ will attract each others, whereas if two state vectors $W_i$ and $W_j$ are dissimilar, then $x_i$ and $x_j$ will repel each other. For definiteness, we set
\begin{align}\label{C-3}
\begin{cases}
\Psi_a^{ij} = \Psi_a(W_{i},W_{j})=(S_i+R_i)(S_j+R_j)+I_i I_j+\varepsilon_a, \quad i, j \in [N], \\
\Psi_r^{ij} = \Psi_r(W_{i},W_{j})=(S_i+R_i)I_j+(S_j+R_j)I_i+\varepsilon_r,  \quad i, j \in [N],
\end{cases}
\end{align}
where $\alpha$ and $\beta$ are positive constants. Here $\varepsilon_a$ and $\varepsilon_r$ are positive constants presenting social distancing. These social distancing constants play an important role in preventing collisions between particles.  Note that the ansatz \eqref{C-3} satisfies symmetry:
\begin{equation} \label{C-3-1}
\Psi_a^{ij} = \Psi_a^{ji}, \quad \Psi_r^{ij} = \Psi_r^{ji}, \quad i, j \in [N]. 
\end{equation}
Next, we will impose conditions on $\alpha$ and $\beta$ later in \eqref{C-2}. To illustrate the functional relations in the right-hand side of \eqref{C-3}, we consider an ensemble which is partitioned  into two sub-ensembles:
 \[ \{ \mbox{Infected particles} \} \quad \mbox{or} \quad  \{ \mbox{susceptible or recovered particles} \}. \]
Recall that $I_i$ and $S_i + R_i$ are probabilities that $i$-th particle are in infected state and are not in infected state, respectively. \newline

Consider the inner product-like function $G$  as follows:
\begin{align}
\begin{aligned} \label{C-4}
G(W_i, W_j) &= (S_i + R_i, I_i) \cdot (S_j + R_j, I_j) = (S_i+R_i)(S_j+R_j)+I_i I_j \\
&=(1-I_i)(1-I_j)+I_iI_j.
\end{aligned}
\end{align}
Note that $G(W_i, W_j)$ becomes larger when $I_i\simeq I_j$ and $S_i+R_i\simeq S_j+R_j$. If $W_i$ and $W_j$ are similar, then $G(W_i, W_j)$ becomes larger, and one has
\begin{equation} \label{C-5}
0\leq G(W_i, W_j)=(1-I_i)(1-I_j)+I_iI_j\leq \big[(1-I_i)+I_i\big]\cdot\big[(1-I_j)+I_i\big]=1.
\end{equation}
Now, we set
\begin{equation} \label{C-6}
\Psi_a^{ij} =\varepsilon_a + G(W_i, W_j) \geq 0, \quad 
\Psi_r^{ij} =1+\varepsilon_r -G(W_i, W_j) \geq 0. 
\end{equation}
Then, $\Psi_a$ and $\Psi_r$ satisfy the following monotonicity properties: 
\begin{enumerate}
\item
If $W_i$ and $W_j$ become similar,  $\Psi_a$ increases and $\Psi_r$ decreases.

\vspace{0.1cm}

\item 
If $W_i$ and $W_j$ become dissimilar, $\Psi_a$ decreases and $\Psi_r$ increases.
\end{enumerate}

\vspace{0.2cm}

\noindent The relations \eqref{C-4} and \eqref{C-6} yield \eqref{C-3}. On the other hand, it follows from \eqref{C-3} and  \eqref{C-5} that the coupling weight functions $\Psi_a$ and $\Psi_r$ admit positive lower bound and upper bounds:
\begin{equation} \label{C-6-0}
0<\varepsilon_a \leq \Psi_a \leq 1+\varepsilon_a =:M_a, \qquad 
0<\varepsilon_r \leq \Psi_r \leq 1+\varepsilon_r =:M_r.
\end{equation}
\begin{proposition}
Suppose that system parameters satisfy 
\begin{equation} \label{C-7}
 \Psi=0,\quad \kappa_1=\frac{a}{N-1},\quad b^i=b, \quad  i \in [N],
 \end{equation}
 and let $(S, I, R)$ and $(S_i, I_i, R_i)$ be the solutions of system \eqref{B-1} and system \eqref{A-1} with the initial data:
\begin{align*}
\begin{aligned}
& S(0)=S^0,\quad I(0)=I^0,\quad R(0)=R^0, \\
&  S_i^0=S^0, \quad I_i^0=I^0,\quad R_i^0=R^0, \quad i \in [N].
\end{aligned}
\end{align*}
Then we have 
\begin{equation} \label{C-8}
S_i(t)=S(t),\quad I_i(t)=I(t),\quad R_i(t)=R(t), \quad t \geq 0,~~ i \in [N].
\end{equation}
\end{proposition}
\begin{proof} Suppose the conditions \eqref{C-7} hold. Then, the relation \eqref{C-1} and \eqref{B-1}  become 
\[ a^{ij}=\begin{cases}
1 \quad&\text{if}\quad i\neq j,\\
0&\text{otherwise},
\end{cases}
\quad \mbox{and} \quad 
\begin{cases} 
\displaystyle \dot{S}_i =-\frac{a}{N-1}   \sum_{\substack{1\leq j\leq N\\j\neq i}}  S_i I_j,  \quad t > 0, \\
\displaystyle \dot{I}_i = \frac{a}{N-1}  \sum_{\substack{1\leq j\leq N\\j\neq i}}  S_i I_j-b I_i, \\
\displaystyle \dot{R}_i =b I_i, \quad i \in [N].
\end{cases}
\]
Since $(S_i, I_i, R_i)$ and $(S, I, R)$ satisfy the same system  \eqref{B-1}  with the same initial data, by the uniqueness of ODE solution, we have the desired uniqueness \eqref{C-8}. 
\end{proof}
%
%
%

\subsection{Symptom expression vector} \label{sec:3.3}
Next, we introduce the symptom expression vector $s=(s_1, \cdots, s_N)$ motivated by the ongoing pandemic COVID-19. In April 2020, the daily COVID-19 confirmed number per day of Korea was less than 20. Most of confirmed people were isolated, however patients with no symptoms of COVID-19 were still spreading the virus. That is why the daily confirmed number per day does not goes to 0 directly. So we need to consider the ratio of symptom expressions. Some people show symptoms well when they got the virus, in contrast some people does not show any symptoms although they got the virus. So even for the same inputs, different degree of output can emerge. Thus, we define the symptom expression vector that are related to previous phenomena. We define define each component $s_i \in [0, 1]$ to following the following properties: \newline
\begin{itemize}
\item If the $i$-th particle shows the symptom well, then we put large value to $s_i$.

\vspace{0.1cm}

\item If  the $i$-th particle does not show the symptom well, then we put small value to $s_i$.
\end{itemize}

\vspace{0.1cm}

We may use the symptom expression vector $s$ to express the condition of being suspected as a patient. In this paper, we assume that the symptom expression vector $s$ is a time-independent constant vector. We assume that the condition of being suspected as a patient only depends on the product $s_iI_i$, and we also define some fixed threshold constant $c$ to make decision. If $s_iI_i(t)\geq c$, then the $i^{th}$ person will be confirmed at time $t$.
\section{Relaxation of epidemic states}\label{sec:4}
\setcounter{equation}{0}
In this section, we study the relaxation dynamics of the SIR-flock model \eqref{A-1}. First, we show that $\mathcal{S}$ is invariant along \eqref{A-1}. 
\begin{lemma} \label{L4.1}
Let $\{(W_i, x_i)\}$ be a solution of system \eqref{A-1}. Then, one has 
\begin{equation} \label{NNN-0}
W_i(t)\in\mathcal{S}, \qquad  \frac{d}{dt} \Big( \sum_{i=1}^N{x_i} \Big) =0, \quad t > 0, \quad i \in [N],
\end{equation}
where ${\mathcal S}$ is the set defined in \eqref{A-0-1}. 
\end{lemma}
\begin{proof} We basically use the same arguments as in Proposition \ref{P2.1}. \newline

\noindent (i)~(Verification of the second relation in \eqref{NNN-0}):~It follows from $\eqref{A-1}_1$ that 
\[
\frac{d}{dt} (S_i + I_i + R_i) =0, \quad t > 0, \quad i \in [N].
\]
This yields
\begin{equation} \label{NNN-1}
S_i(t)+I_i(t)+R_i(t)=1, \quad t \geq 0.
\end{equation}

\vspace{0.2cm}

Now, we need to show the positivity of $S_i, I_i$ and $R_i$.  \newline

\noindent $\bullet$~(Positivity of $S_i$):~Due to \eqref{NNN-1}, it suffices to check the positivity of $S_i, I_i$ and $R_i$. Now, we integrate $\eqref{A-1}_1$ to find 
\[
S_i(t)=S_i^0 \exp\left(-\kappa_1 \int_0^t \sum_{j=1}^N a^{ij}(\tau)I_j(\tau)d\tau\right) \geq 0, \quad \mbox{for}~ t \geq 0. 
\]
Thus, one has 
\[
S_i(t)\geq0, \qquad \mbox{for}~t \geq 0, \quad i \in [N].
\]

\vspace{0.2cm}

\noindent $\bullet$~(Positivity of $I_i$):~It follows from \eqref{B-1}$_2$ that 
\[
\dot{I}_i+b^i I_i= \kappa_1 \sum_{j=1}^N a^{ij}S_iI_j.
\]
This yields
\[
\frac{d}{dt}(I_ie^{b^i t})= \kappa_1 \sum_{j=1}^N a^{ij}S_iI_j e^{b^i t}.
\]
Now, we introduce a time-varying minimal index $m_t \in [N]$ such that 
\[
I_{m_t}(t):=\min_{1 \leq i \leq N}I_i(t).
\]
Since each $I_i(t)$ is analytic, there exists the refinement of time-interval 
\[ 0=t_0<t_1<t_2<\cdots \]
such that the index of $I_{m_t}$ is not changed on each interval $[t_j, t_{j+1})$, and 
\[
\frac{d}{dt}(I_{m_t} e^{b^{m_t} t})= \kappa_1 \sum_{j=1}^N a^{m_t j}S_{m_t} I_j e^{b^{m_t} t}, \quad t > 0. 
\]
Therefore, we have
\begin{align*}
\frac{d}{dt}(I_{m_t}e^{b^{m_t} t}) = \kappa_1 \sum_{j=1}^N a^{m_t j}S_{m_t}I_j e^{b^{m_t} t}
\geq \kappa_1 \sum_{j=1}^N a^{m_tj}S_{m_t}I_{m_t} e^{b^{m_t} t}
= I_{m_t} e^{b^{m_t} t} \Big( \kappa_1 \sum_{j=1}^N a^{m_t j}S_{m_t} \Big).    
\end{align*}
This yields
\begin{align*}
I_{m_t}e^{b^{m_t} t} \geq I_{m_0} \exp\left( \kappa_1 \sum_{j=1}^N \int_0^t \ a^{m_\tau j}(\tau)S_{m_\tau} d\tau\right)\geq 0.
\end{align*}
Therefore, one has 
\[
I_i(t) \geq I_{m_t} (t) \geq 0, \quad \forall~i \in [N].
\]
\vspace{0.5cm}

\noindent $\bullet$~(Positivity of $R_i$): Since $S_i + I_i + R_i = 1,~S_i \geq 0$ and $I_i \geq 0$, one has 
\[ R_i \geq 0. \]
Therefore, we combine all the estimates for $S_i, I_i$ and $R_i$ to get 
\[ W_i \in {\mathcal S}. \]

\vspace{0.5cm}

\noindent (ii)~(Verification of the second relation in \eqref{NNN-0}):~Recall the equation for $x_i$: 
\[ 
 \dot{x}_i =\displaystyle\frac{\kappa_2}{N}  \sum_{j\neq i} \left(\Psi_a^{ij}\frac{x_j-x_i}{\|x_j-x_i\|^\alpha}\right)-\frac{\kappa_3}{N}  \sum_{j\neq i}\left(\Psi_r^{ij}\frac{x_j-x_i}{\|x_j-x_i \|^\beta}\right). \]
 Since the right-hand side of the above equation is skew-symmetry with respect to interchange map $i~\longleftrightarrow~j$ using \eqref{C-3-1},  the total sum $\sum_{i} x_i$ is time-invariant. 
\end{proof}
Now, we discuss the asymptotic convergence of the probability vector $W_i$ toward a constant probability vector $W^{\infty} \in {\mathcal S}$ in \eqref{A-0-1}. 
\begin{theorem} \label{T4.1}
Let $\{ (W_i,x_i) \}$ be a solution of system \eqref{A-1}. Then, the following assertions hold.
\begin{enumerate}
\item
There exists a constant state $\{ (S_i^{\infty}, 0, R_i^{\infty})\}$ with $0\leq S^{\infty}_i,~R^{\infty}_i\leq 1$ such that 
\[
\lim_{t\rightarrow \infty} (S_i(t), I_i(t), R_i(t)) = (S^{\infty}_i, 0, R^{\infty}_i), \quad i \in [N].
\]
\item
If the recovering value $b_i$ satisfies  
\begin{equation} \label{New-1}
 \min_{1 \leq i \leq N} b^i > \frac{\kappa_1(N-1)}{L^\gamma},
 \end{equation}
then there exists a positive constant $\displaystyle \lambda := \min_{1 \leq i \leq N} b^i  -\frac{\kappa_1(N-1)}{L^\gamma} >0$ such that 
\[
\Big| \frac{1}{N} \sum_{i=1}^N{I_i(t)} \Big|\leq e^{-\lambda t}.
\]
\end{enumerate}
\end{theorem}
\begin{proof}
\noindent (i)~Since $I_j$ and $S_i$ are non-negative,
\[
\dot{S}_i=-\kappa_1 \left(\sum_{j=1}^N a^{ij}I_j\right)S_i \leq 0,
\]
which means $S_i$ is non-increasing. Since $S_i$ is bounded below by zero, there are some constants $S^{\infty}_i$ that $S_i(t)$ converges to $S_i^{\infty}$, as $t$ goes $\infty$ by monotone convergence theorem. Similarly,
\[
\dot{R}_i=b^iI_i \geq 0,
\]
and $R_i$ is non-decreasing and bounded above, there are some constants $R^{\infty}_i$ that $R_i(t)$ converges to $R_i^{\infty}$, as $t$ goes $\infty$. Moreover, as $R_i$ converges,
\begin{align*}
0=\lim_{t\rightarrow \infty}\dot{R}_i(t)=\lim_{t\rightarrow \infty}b^i I_i(t),\quad\text{thus}\quad\lim_{t\rightarrow \infty}I_i(t)=0.
\end{align*}

\vspace{0.2cm}

\noindent (ii)~Recall the equation for $I_i$:
\[ \dot{I}_i + b^i I_i = \kappa_1 \sum_{j=1}^N a^{ij}S_i I_j. \]
Now, we use the method of integrating factor to derive Gronwall's inequality for $\sum_{i=1}^N{e^{b^i t}I_i}$:
\[
\frac{d}{dt}\left(\sum_{i=1}^{N}{e^{b^i t}I_i}\right) 
= \kappa_1 e^{b^i t}\sum_{i,j\in\mathcal{N}}{a^{ij}S_iI_j} 
\leq e^{b^it}\sum_{\substack{i,j\in\mathcal{N}\\ j\neq i}}{\frac{\kappa_1}{L^\gamma}I_j} 
= \frac{\kappa_1(N-1)}{L^\gamma}\left(\sum_{i=1}^{N}{e^{b^i t}I_i}\right).
\]
Thus, we have
\begin{align*}
\left(\sum_{i=1}^{N}{e^{b^i t}I_i}\right) \leq \left(\sum_{i=1}^{N}{I_i(0)}\right)e^{\frac{\kappa_1(N-1)}{L^\gamma} t}\leq Ne^{\frac{\kappa_1(N-1)}{L^\gamma} t}.
\end{align*}
This implies
\begin{align*}
I_i \leq \sum_{i=1}^{N}{I_i} \leq Ne^{\left(\frac{\kappa_1(N-1)}{L^\gamma}-b^i \right)t} \leq Ne^{\left(\frac{\kappa_1(N-1)}{L^\gamma}-\min_{1 \leq i \leq N} b^i \right)t} =:Ne^{-\lambda t}.
\end{align*}
\end{proof}
Next, under a more relaxed condition compared to \eqref{New-1}, we improve the second result of Theorem \ref{T4.1} as follows.
\begin{corollary}
Let $(W_i, x_i)$ be the solution of the system \eqref{A-1}. If $\{ b^i \}$ satisfies following  condition:
\begin{equation} \label{New-2}
 \min_{1 \leq i \leq N} b^i > \frac{\kappa_1}{L^\gamma} \Big( \max_{1 \leq i \leq N} \sum_{j\neq i}S_j^0 \Big),
\end{equation}
then $\sum_{i=1}^N{I_i}$ decays to zero exponentially fast.
\end{corollary}
\begin{proof}
First, we write the dynamics $\eqref{A-1}_1$ of $I_i$s in matrix form:
\begin{align*}
    \dot{\mathbf{I}}^\top &=\kappa_1\text{diag}(S_1,S_2,\cdots,S_N)\mathbf{A}\mathbf{I}^\top -b\mathbf{I}^\top ,
\end{align*}
where $\mathbf{I}$ and $ \mathbf{A}$ are given as follows.
\[ \mathbf{I}(t):=(I_1(t),\cdots,I_N(t)) \quad \mbox{and} \quad \mathbf{A}(t):=(a^{ij}(t))_{1\leq i, j\leq N}. \] 
Since each of $S_i$ is decreasing, 
\[ S_i(t)\leq S_i(0) \quad \mbox{and} \quad a^{ij}\leq L^{-\gamma} \quad \mbox{for $i\neq j$}. \]
This implies that for every $t$,
\[
\text{diag}(S_1(t),S_2(t),\cdots,S_N(t))\mathbf{A}(t)\leq \mathbf{\bar{A}},
\]
where $\leq$ is a partial order compoentwise, and 
\begin{align*}
\mathbf{\bar{A}}_{ij}=
    \begin{cases}
    S_i^0 L^{-\gamma}, & i\neq j,\\
    0, & i=j.
    \end{cases}
\end{align*}
Therfore, we obtain
\begin{equation} \label{New-3}
    \dot{\mathbf{I}}^\top  =\kappa_1\text{diag}(S_1,S_2,\cdots,S_N)\mathbf{A}\mathbf{I}^\top -b\mathbf{I}^\top 
    \leq (\kappa_1 \mathbf{\bar{A}}-bI)\mathbf{I}^\top .
\end{equation}
We multiply the both sides of \eqref{New-3} by $\mathbbm{1}$ to yield
\begin{align*}
    \frac{d}{dt}\left(\sum_{i=1}^N{I_i(t)}\right)\leq \sum_{i=1}^N \lambda_i I_i(t), 
\end{align*}
where $\lambda_i$ is given as follows.
\[  \lambda_i :=\kappa_1 L^{-\gamma} \sum_{j\neq i}S_j^0 - b^i. \]
If 
 \[ \displaystyle \min_{1 \leq i \leq N} b^i > \kappa_1 L^{-\gamma} \sum_{j\neq i}S_j^0, \]
 then $\lambda_i$ is negative, so $\sum_{i=1}^N{I_i(t)}$ decays to zero exponentially fast.
\end{proof}
\begin{remark}
Since 
\[ N-1 \geq  \max_{1 \leq i \leq N} \sum_{j\neq i}S_j^0, \]
the condition \eqref{New-2} gives a more relaxed condition on $\{ b^i \}$ compared to \eqref{New-1}. 
\end{remark}
\section{Quantitative estimates for relative distances} \label{sec:5}
\setcounter{equation}{0}
In this section, we provide explicit quantitative estimates on relative distances between particles. For a given spatial configuration $\{ x_i \}$, we consider the set of all relative distances $\{\|x_i - x_j \| \}$ for $i \not = j$, and rearrange them in an increasing order:
\[  D_1(t)\leq D_2(t)\leq \cdots \leq D_{\mathcal Q}(t) =: D(t), \quad \mbox{where}~~{\mathcal Q} :=  \frac{N(N-1)}{2}. \]
In what follows, we derive a uniform lower-bound for $D_1$ and a uniform upper-bound for $D_{\mathcal Q}$ so that we can control the singular terms in \eqref{A-1} uniformly in time. First, we recall notation:
\begin{equation} \label{B-7}
\Lambda := \Big( \frac{1}{2} \Big)^{{\mathcal Q}} \max \{  D_1(0), 1 \}^{\mathcal Q} \Big(  \Big( \max\Big\{ \Big( \frac{m_r}{ 2 M_a} \Big)^{\frac{1}{\beta -\alpha}},   1\Big \} \Big)^{\mathcal Q} \Big( \frac{m_r}{2N (M_a + M_r)} \Big )^{\frac{{\mathcal Q}}{\beta - 1}}.
\end{equation}

\begin{theorem} \label{T5.1}
Let $\{ (W_i,x_i) \}$ be a solution of system \eqref{A-1} with the initial data $\{ (W^0_i,x^0_i) \}$. Then, the following assertions hold.
\begin{enumerate}
\item
(Existence of a positive lower bound to $D_1$):~If initial data satisfy 
\[ \min_{1 \leq i \not = j \leq N} \|x_i^0 -x_j^0 \| > 0, \]
there is a positive constant $\delta_1 >0$ such that
\[
\inf_{0\leq t<\infty} D_1(t) \geq\delta_1.
\]
\item
(Existence of a positive upper bound to $D_{\mathcal Q}$):
If initial data and system functions satisfy
\[ \Lambda > \left(\frac{\kappa_3 M_r}{\kappa_2 m_a}\right)^{\beta-\alpha},
\]
there exists a positive constant $\delta_\infty$ such that
\[
\sup_{0\leq t<\infty} D_{\mathcal Q}(t) \leq \delta_\infty.
\]
\end{enumerate}
\end{theorem}
\begin{proof}
We leave its proof in the next two subsections. 
\end{proof}

\subsection{A positive lower bound for minimal relative distance} \label{sec:5.1}
In this subsection, we show that $D_1(t)$ has a positive lower bound $\delta_1$ which is independent of $t$. Since the particles do not collide in finite-time interval, there exists an analytic solution in finite time (see Theorem \ref{T2.2}). Therefore, each distance $D_{ij}$ and ordered distance $D_i$ are Lipschitz continuous. In addition, according to the analyticity,  for given $D_{ij}$, there exist a sequence of times $(t_n)$: 
\[
0=:t_0<t_1<t_2<\cdots <t_n<\cdots,
\]
such that we can decompose the whole time interval $[0,\infty)$ into the union of subintervals
\[
[0,\infty)=\bigcup_{l=1}^{\infty}{T_l},\quad T_l=[t_{l-1},t_l),\quad l\geq 1.
\]
to make sure that $D_{ij}$ is not changed in each subinterval. To demonstrate the existence of positive uniform lower bound, we first provide three lemmas. First, we show a positive lower bound of maximal distance $D_{\mathcal Q}(t)$:
\begin{align}\label{EQ4.1}
D(t):=D_{\mathcal Q}(t)=\max_{1 \leq i \not = j \leq N} \|x_i(t) - x_j(t) \|.    
\end{align}
Next, we provide a series of lemmas to estimate maximal and minimal relative distances. First, we derive a differential inequality for $D$.
\begin{lemma}\label{L5.2}
Let $\{ (W_i,x_i) \}$ be a solution to system \eqref{A-1}. Then, the functional $D$ in \eqref{EQ4.1} satisfies 
\begin{align*}
\frac{d}{dt}{D(t)}^2&\geq
-\frac{4\kappa_2 M_a}{N\cdot D^{\alpha-2}}+\frac{4\kappa_3 m_r}{N\cdot D^{\beta-2}}
+\frac{2}{N}\sum_{k\neq i,j}\left[\frac{(x_i-x_j)\cdot(x_k-x_i)}{\|x_k-x_i \|^\alpha}\left(\kappa_2 M_a - \frac{\kappa_3 m_r}{\|x_k-x_i\|^{\beta-\alpha}}\right)\right]\\
&+\frac{2}{N}\sum_{k\neq i,j}\left[\frac{(x_i-x_j)\cdot(x_j-x_k)}{\|x_k-x_j \|^\alpha}\left(\kappa_2 M_a - \frac{\kappa_3 m_r}{ \|x_k-x_j \|^{\beta-\alpha}}\right)\right].
\end{align*}
\end{lemma}
\begin{proof}
\noindent For given $l\in \mathbb{N}\cup\{0\}$, we choose $i$ and $j$ such that 
\[
D(t)= \|x_i(t)-x_j(t) \|\quad \forall~ t\in T_l.
\]
Then, by direct calculation, we obtain
\begin{align}
\begin{aligned} \label{E-4-1}
\frac{d}{dt}{D(t)}^2 &=\frac{d}{dt} \|x_i-x_j \|^2=2(x_i-x_j)\cdot(\dot{x_i}-\dot{x_j})\\
&=2(x_i-x_j)\cdot \left[\frac{\kappa_2}{N}\sum_{k\neq i}{\left(\Psi_a^{ik} \frac{x_k-x_i}{\|x_k-x_i \|^\alpha}\right)-\frac{\kappa_3}{N}\sum_{k\neq i}\left(\Psi_r^{ik}\frac{x_k-x_i}{\|x_k-x_i\|^\beta}\right)}\right]\\
&\hspace{0.2cm}-2(x_i-x_j)\cdot \left[\frac{\kappa_2}{N}\sum_{k\neq j}{\left(\Psi_a^{jk}\frac{x_k-x_j}{\|x_k-x_j \|^\alpha}\right)-\frac{\kappa_3}{N}\sum_{k\neq j}\left(\Psi_r^{jk}\frac{x_k-x_j}{\|x_k-x_j \|^\beta}\right)}\right]\\
&= -\frac{4\kappa_2\Psi^{jk}_a}{N \|x_i-x_j \|^{\alpha-2}}+\frac{4\kappa_3\Psi^{jk}_r}{N \| x_i-x_j \|^{\beta-2}}\\
&\hspace{0.2cm}+\frac{2}{N}\sum_{k\neq i,j}\left[\kappa_2\Psi_a^{ik}\frac{(x_i-x_j)\cdot(x_k-x_i)}{\|x_k-x_i \|^\alpha}-\kappa_3\Psi_r^{ik}\frac{(x_i-x_j)\cdot(x_k-x_i)}{\|x_k-x_i\|^\beta}\right]\\
&\hspace{0.2cm}-\frac{2}{N}\sum_{k\neq i,j}\left[\kappa_2\Psi_a^{jk}\frac{(x_i-x_j)\cdot(x_k-x_j)}{\|x_k-x_j\|^\alpha}-\kappa_3\Psi_r^{jk}\frac{(x_i-x_j)\cdot(x_k-x_j)}{\|x_k-x_j\|^\beta}\right]\\
&=: \mathcal{I}_{21}+\mathcal{I}_{22}+\mathcal{I}_{23},
\end{aligned}
\end{align}
for $t\in T_l$. We now estimate $\mathcal{I}_{2i}$ ($i=1,2,3$) one by one.

\vspace{0.5cm}

\noindent$\bullet$ (Estimate of $\mathcal{I}_{21}$):~We use 
\[ \Psi_a\leq M_a, \quad \Psi_r\geq m_r \quad \mbox{and} \quad \|x_i-x_j\|=D \]
to get 
\begin{equation} \label{E-4-2}
\mathcal{I}_{21}\geq -\frac{4\kappa_2 M_a}{N\cdot D^{\alpha-2}}+\frac{4\kappa_3 m_r}{N\cdot D^{\beta-2}}.
\end{equation}

\vspace{0.2cm}

\noindent$\bullet$ (Estimate of $\mathcal{I}_{22}$) : Since $\|x_i-x_j\|$ is the maximal distance,
\begin{align*}
\begin{aligned} 
(x_i-x_j)\cdot(x_k-x_i) &=(x_i-x_j)\cdot(x_k-x_j+x_j-x_i) \\
&=(x_i-x_j)\cdot(x_k-x_j)-|x_i-x_j|^2\leq 0.
\end{aligned}
\end{align*}
Similarly, we have
\[
(x_i-x_j)\cdot(x_k-x_j)\geq 0.
\]
Therefore, one has
\begin{align}
\begin{aligned} \label{E-4-3}
\mathcal{I}_{22} &\geq \frac{2}{N}\sum_{k\neq i,j}\left[\kappa_2 M_a\frac{(x_i-x_j)\cdot(x_k-x_i)}{\|x_k-x_i\|^\alpha}-\kappa_3 m_r\frac{(x_i-x_j)\cdot(x_k-x_i)}{\|x_k-x_i\|^\beta}\right]\\
&=\frac{2}{N}\sum_{k\neq i,j}\left[\frac{(x_i-x_j)\cdot(x_k-x_i)}{\|x_k-x_i \|^\alpha}\left(\kappa_2 M_a - \frac{\kappa_3 m_r}{\|x_k-x_i \|^{\beta-\alpha}}\right)\right].
\end{aligned}
\end{align}
\vspace{0.5cm}

\noindent$\bullet$ (Estimate of $\mathcal{I}_{23}$):~Similar to the estimate of ${\mathcal I}_{22}$, we obtain
\begin{align}
\begin{aligned} \label{E-4-4}
\mathcal{I}_{23}&\geq \frac{2}{N}\sum_{k\neq i,j}\left[\kappa_2 M_a\frac{(x_i-x_j)\cdot(x_j-x_k)}{\|x_k-x_j \|^\alpha}-\kappa_3 m_r\frac{(x_i-x_j)\cdot(x_j-x_k)}{\|x_k-x_j\|^\beta}\right]\\
&=\frac{2}{N}\sum_{k\neq i,j}\left[\frac{(x_i-x_j)\cdot(x_j-x_k)}{\|x_k-x_j \|^\alpha}\left(\kappa_2 M_a - \frac{\kappa_3 m_r}{\|x_k-x_j \|^{\beta-\alpha}}\right)\right].
\end{aligned}
\end{align}
In \eqref{E-4-1}, we combine all the estimates \eqref{E-4-2}, \eqref{E-4-3} and \eqref{E-4-4} to find the desired estimate. 
\end{proof}
To derive a positive lower bound for $D_1$, we first verify that the maximal distance $D$ is bounded away from zero uniformly in time, and then we show that this positive lower bound propagates to the lower graded relative distance when we reach to the minimal relative distance $D_1$ (the first assertion in Theorem \ref{T5.1}) via the following two lemmas. \newline
Now we derive a positive lower bound for the maximal distance $D(t)$.
\begin{lemma}\label{L5.2}
\emph{(maximal relative distance)}
Suppose that initial spatial configuration is non-collisional:
\[ \min_{1 \leq i \not = j \leq N}  \|x_i^0 -x_j^0 \| > 0, \]
and let $\{(W_i,x_i)\}$ be a solution to system \eqref{A-1}. Then, one has 
\[
\inf_{0\leq t<\infty} D(t) \geq \delta_{\mathcal Q} := \min\left\{D(0),~~\left(\frac{\kappa_3 m_r}{\kappa_2 M_a}\right)^{\frac{1}{\beta-\alpha}}\right\}.
\]
\end{lemma}
\begin{proof} 
We use an induction argument on the time intervals $T_l=[t_{l-1},t_l)$. \newline

\noindent$\bullet$ (Initial step): we claim that 
\[
D(t)\geq\delta_{\mathcal Q},\quad \text{for } t\in T_1=[0,t_1).
\]
Suppose not, since $\delta_{\mathcal Q} \leq D(0)$, there exist $t_{11}$ and $t_{12}$ such that 
\[ 0\leq t_{11}<t_{12}<t_1, \quad  D(t_{11})=\delta_{\mathcal Q} \quad  \text{and}\quad D(t)<\delta_{\mathcal Q} \text{ in } (t_{11},t_{12}].
\]
Therefore, for $t\in(t_{11},t_{12}]$,
\begin{align*}
D<\delta_{\mathcal Q} \leq \left(\frac{\kappa_3 m_r}{\kappa_2 M_a}\right)^{\frac{1}{\beta-\alpha}}.
\end{align*}
For all $k\neq i,j$, we get
\[
-\frac{4\kappa_2 M_a}{N\cdot D^{\alpha-2}}+\frac{4\kappa_3 m_r}{N\cdot D^{\beta-2}} \geq 0, \quad \kappa_2 M_a - \frac{\kappa_3 m_r}{|x_k-x_i|^{\beta-\alpha}}\leq0, \quad \kappa_2 M_a - \frac{\kappa_3 m_r}{|x_k-x_j|^{\beta-\alpha}}\leq0,
\]
By Lemma \ref{L4.1}, one has 
\[
\frac{d}{dt}D^2 \geq 0\quad\text{for }t\in(t_{11},t_{12}].
\]
This implies
\[
D(t_{12})\geq D(t_{11})=\delta_{\mathcal Q},
\]
which is contradictory to 
\[
D(t)<\delta_{\mathcal Q} \quad\text{in}\quad (t_{11},t_{12}].
\]
\vspace{0.2cm}

\noindent$\bullet$ (Inductive step) : Suppose that for some $l\geq1$,
\[
D(t)\geq\delta_{\mathcal Q} \quad\text{for}\quad t\in T_l=[t_{l-1},t_l).
\]
Now we consider $D(t)$ on $T_{l+1}$. Due to the continuity of $D(t)$ and induction hypothesis, one has
\[
D(t_l)=\lim_{t\rightarrow t_l^-}{D(t)}\geq\delta_{\mathcal Q}.
\]
Therefore, we can use the same criteria above to get
\[
D(t)\geq\delta_{\mathcal Q} \quad\text{for}\quad t\in T_{l+1}.
\]
Thus, one has the desired estimate.
\end{proof}

\vspace{0.2cm}

In the following lemma, we show the backward propagation of a positive lower bound for the ordered relative distance. 
\begin{lemma}\label{L5.3}
\emph{(Backward propagation of lower bounds)}
Suppose that initial spatial configuration is non-collisional:
\[  \min_{1 \leq i \not = j \leq N}  \|x_i^0 -x_j^0 \| > 0, \]
and let $\{(W_i,x_i)\}$ be a solution to system \eqref{A-1}, and  let $ p \in [1, {\mathcal Q})$ be a fixed constant. If there exist a positive constant $\delta_q$  such that
\[
\inf_{0\leq t<\infty}D_q(t)\geq \delta_q \text{ for all }~ q \in (p, {\mathcal Q}],
\]
then there exists a positive constant $\delta_p$ such that 
\[
\inf_{0\leq t<\infty}D_p(t) \geq \delta_p.
\]
\end{lemma}
\begin{proof}
We refer to the proof of  Lemma 3.4 in \cite{HJKPZ}.
\end{proof}

\vspace{0.5cm}

Now, we are ready to provide a positive lower bound for the minimal relative distance $D_1$.  \newline

\noindent {\bf Proof of the first assertion in Theorem \ref{T5.1}}:~Suppose initial spatial configuration is noncollisional:
\[ \min_{1 \leq i \not = j \leq N}  \|x_i^0 -x_j^0 \| > 0, \]
and let $\{(W_i,x_i)\}$ be a solution to system \eqref{A-1}. By Lemma \ref{L5.2}, there exists $\delta_{\mathcal Q} > 0$ such that 
\[ \inf_{0\leq t<\infty} D(t) \geq \delta_{\mathcal Q}.   \] 
For $q = {\mathcal Q}$, we apply Lemma \ref{L5.3} to show that 
\[ \exists~\delta_{{\mathcal Q}-1} > 0 \quad \mbox{such that} \quad \inf_{0\leq t<\infty}D_{{\mathcal Q}-1}(t) \geq \delta_{{\mathcal Q}-1}. \]
Again, we apply Lemma \ref{L5.3} inductively until we reach $q = 1$ to derive the desired a positive lower bound for $D_1$.  $\qed$ 

\subsection{ A positive upper bound for maximal relative distance} \label{sec:5.2} In this subsection, we derive an existence of a positive upper bound for the maximal relative distance. 
First, we study the estimate for a differential inequality. 
\begin{lemma}\label{L5.4}
Let $y: [0,\infty)\rightarrow \mathbb{R}^+$ be a differentiable function satisfying the following differential inequality:
\begin{align*}
\begin{cases}
y'\leq -ay^{-p}+by^{-q}, \quad t > 0,\\
y(0)=y^0,
\end{cases}
\end{align*}
where constants $a, b, p$ and $q$ satisfy
\[
a > 0, \quad b > 0, \quad 0\leq p<q,
\]
Then, there exists a positive constant $\delta$ such that 
\[ \sup_{0 \leq t < \infty} y(t) \leq \delta. \]
\end{lemma}
\begin{proof}
We use phase line analysis for ${\dot z} =-az^{-p}+bz^{-q}$. For this, we define
\[
y^*:=\left(\frac{b}{a}\right)^{\frac{1}{q-p}}.
\]
Then, it is easy to see that 
\[
f(y)>0 \quad  \text{ for } y<y^*, \quad  f(y^*)=0 \quad \text{ and } \quad  f(y)<0 \quad \text{ for } y>y^*.
\]
Therefore, $y$ is uniformly bounded by $\max\{y^0,y^*\}$.
\end{proof}
Now, we are ready to provide a proof for the second assertion in Theorem \ref{T5.1}. \newline

\noindent {\bf Proof of the second assertion in Theorem \ref{T5.1}}: Suppose initial data and system parameters satisfy 
\[ D(0)<\infty,\quad \delta_1 > \left(\frac{\kappa_3 M_r}{\kappa_2 m_a}\right)^{\beta-\alpha}, \]
and  let $\{(W_i,x_i)\}$ be a solution to system \eqref{A-1}. Then,  for $t \in (0, \infty)$, we choose $i = i_t$ and $j = j_t$ such that 
\[ D(t) = \|x_i(t) - x_j(t) \|. \]
By Lemma \ref{L5.2}, one has 
\begin{align}
\begin{aligned} \label{E-5}
& \frac{d}{dt}{D(t)}^2 =\frac{d}{dt} \|x_i-x_j \|^2=2(x_i-x_j)\cdot(\dot{x_i}-\dot{x_j})\\
& \hspace{1cm} =2(x_i-x_j)\cdot \left[\frac{\kappa_2}{N}\sum_{k\neq i}{\left(\Psi_a^{ik}\frac{x_k-x_i}{\|x_k-x_i\|^\alpha}\right)-\frac{\kappa_3}{N}\sum_{k\neq i}\left(\Psi_r^{ik}\frac{x_k-x_i}{\|x_k-x_i\|^\beta}\right)}\right]\\
& \hspace{1.2cm}  -2(x_i-x_j)\cdot \left[\frac{\kappa_2}{N}\sum_{k\neq j}{\left(\Psi_a^{jk}\frac{x_k-x_j}{\|x_k-x_j\|^\alpha}\right)-\frac{\kappa_3}{N}\sum_{k\neq j}\left(\Psi_r^{jk}\frac{x_k-x_j}{\|x_k-x_j \|^\beta}\right)}\right]\\
& \hspace{1cm} = -\frac{4\kappa_2\Psi_a^{jk}}{N \|x_i-x_j \|^{\alpha-2}}+\frac{4\kappa_3\Psi_r^{jk}}{N \|x_i-x_j \|^{\beta-2}}\\
& \hspace{1.2cm}  +\frac{2}{N}\sum_{k\neq i,j}\left[\kappa_2\Psi_a^{ik}\frac{(x_i-x_j)\cdot(x_k-x_i)}{ \|x_k-x_i \|^\alpha}-\kappa_3\Psi_r^{ik}\frac{(x_i-x_j)\cdot(x_k-x_i)}{\|x_k-x_i\|^\beta}\right]\\
& \hspace{1.2cm} -\frac{2}{N}\sum_{k\neq i,j}\left[\kappa_2\Psi_a^{jk}\frac{(x_i-x_j)\cdot(x_k-x_j)}{\|x_k-x_j \|^\alpha}-\kappa_3\Psi_r^{jk}\frac{(x_i-x_j)\cdot(x_k-x_j)}{\|x_k-x_j \|^\beta}\right]\\
& \hspace{1cm} =: \mathcal{I}_{31}+\mathcal{I}_{32}+\mathcal{I}_{33}, \qquad t\in T_l.
\end{aligned}
\end{align}
 In what follows, we estimate $\mathcal{I}_{3i}$ ($i=1,2,3$) one by one. 
 
 \vspace{0.5cm}

\noindent$\bullet$ (Estimate of $\mathcal{I}_{31}$):~We use 
\[ D = \|x_i-x_j \|, \quad  \Psi_a\geq m_a \quad \mbox{and} \quad \Psi_r\leq M_r \]
to get 
\begin{equation} \label{E-6}
\mathcal{I}_{31}\leq -\frac{4\kappa_2 m_a}{N\cdot D^{\alpha-2}}+\frac{4\kappa_3 M_r}{N\cdot D^{\beta-2}}.
\end{equation}

\vspace{0.5cm}

\noindent$\bullet$ (Estimate of $\mathcal{I}_{32}$): We use $D(t) = \|x_i-x_j \|$ to see
\begin{align*}
\begin{aligned}
(x_i-x_j)\cdot(x_k-x_i) &=(x_i-x_j)\cdot(x_k-x_j+x_j-x_i) \\
&=(x_i-x_j)\cdot(x_k-x_j)-\|x_i-x_j \|^2\leq 0.
\end{aligned}
\end{align*}
Now, we define the angle between $x_i-x_j$ and $x_k-x_i$ by $\theta_{ijk}$. Since
\[
\|x_i-x_j \| \cdot \|x_k-x_i \|\cos(\theta_{ijk})=(x_i-x_j)\cdot(x_k-x_i)\leq 0,
\]
one has 
\[  \cos(\theta_{ijk}) \leq 0. \] 
Thus, one has 
\begin{align*}
\mathcal{I}_{32}&\leq \frac{2}{N}\sum_{k\neq i,j}\left[\kappa_2 m_a\frac{(x_i-x_j)\cdot(x_k-x_i)}{ \|x_k-x_i \|^\alpha}-\kappa_3 M_r\frac{(x_i-x_j)\cdot(x_k-x_i)}{\|x_k-x_i \|^\beta}\right]\\
&= \frac{2}{N}\sum_{k\neq i,j}\left[\kappa_2 m_a\frac{\|x_i-x_j \| \cdot \|x_k-x_i \|\cos{\theta_{ijk}}}{\|x_k-x_i \|^\alpha}-\kappa_3 M_r\frac{\|x_i-x_j \| \cdot \|x_k-x_i \|\cos{\theta_{ijk}}}{\|x_k-x_i \|^\beta}\right]\\
&= \frac{2}{N}\sum_{k\neq i,j}\left[\kappa_2 m_a\frac{\|x_i-x_j \|\cos{\theta_{ijk}}}{\|x_k-x_i \|^{\alpha-1}}-\kappa_3 M_r\frac{\|x_i-x_j \|\cos{\theta_{ijk}}}{\|x_k-x_i \|^{\beta-1}}\right]\\
&=\frac{2}{N}\sum_{k\neq i,j} \|x_i-x_j \|\cos{\theta_{ijk}}\left[\frac{\kappa_2 m_a}{\|x_k-x_i \|^{\alpha-1}}-\frac{\kappa_3 M_r}{\|x_k-x_i \|^{\beta-1}}\right].
\end{align*}
Now, we define the function
\[
f(d):=\frac{\kappa_2 m_a}{d^{\alpha-1}}-\frac{\kappa_3 M_r}{d^{\beta-1}}.
\]
Then, we have
\[ f(d) > 0, \quad \forall~d > \left(\frac{\kappa_3 M_r}{\kappa_2 m_a}\right)^{\beta-\alpha}. \]
Since 
\[
\left(\frac{\kappa_3 M_r}{\kappa_2 m_a}\right)^{\beta-\alpha}<\delta_1\leq |x_k-x_i|,
\]
we have
\begin{equation} \label{E-7}
\mathcal{I}_{32}\leq 0.
\end{equation}

\vspace{0.5cm}

\noindent$\bullet$ (Estimate of $\mathcal{I}_{33}$): Similar to ${\mathcal I}_{32}$, we have 
\begin{equation} \label{E-8}
 \mathcal{I}_{33}\leq 0. 
 \end{equation}

\vspace{0.2cm}

\noindent Finally, in \eqref{E-5}, we combine all the estimates \eqref{E-6}, \eqref{E-7} and \eqref{E-8} to find 
\[
\frac{d}{dt}D^2 \leq-\frac{4\kappa_2 m_a}{N\cdot D^{\alpha-2}}+\frac{4\kappa_3 M_r}{N\cdot D^{\beta-2}},
\]
or equivalently
\[
\frac{d}{dt}D \leq-\frac{2\kappa_2 m_a}{N\cdot D^{\alpha-1}}+\frac{2\kappa_3 M_r}{N\cdot D^{\beta-1}}.
\]
Now, we apply Lemma \ref{L5.4} with parameters
\[
a=\frac{2\kappa_2 m_a}{N}, \quad b=\frac{2\kappa_3 M_r}{N}, \quad p=\alpha-1, \quad q=\beta-1 
\]
to derive the desired estimate.  $\quad \qed$ 
\begin{remark}
Note that for $1\leq\alpha<\beta$, nonexistence of finite-time collisions and uniform lower bound of diameter are always warranted, however conditions on parameter is needed to guarantee uniform upper bound of diameter.
\end{remark}

\section{Relaxation of spatial configuration} \label{sec:6}
\setcounter{equation}{0}
In this section, we study the relaxation of spatial configuration toward a constant configuration. \newline

Consider a perturbed gradient system with decaying forcing $f$:
\begin{equation} \label{F-1}
\dot{x}(t) =-\nabla_xV(x(t))+f(t),  \quad t > 0,
\end{equation}
where $V$ is a one-body potential. 
\begin{proposition}\label{P4.9}
Suppose that external forcing $f$ and potential $V$ satisfy the following conditions:
\begin{enumerate}
\item
 There exist positive constants $C$ and $\lambda$ such that
\[
|f(t)|\leq Ce^{-\lambda t},
\]
\item 
There exists a compact set $K$ such that $x(t)$ is contained in $K$ for all $t\geq 0$,
\item 
V is analytic in $K$ and  the map $t \mapsto |\nabla_xV(x(t))|^2$ is uniformly continuous in time,
\end{enumerate}
and let $x=x(t)$ be a solution of \eqref{F-1}.  Then, there exists $x^\infty \in K$ such that 
\[
\lim_{t\rightarrow \infty}x(t)=x^\infty.
\]
\end{proposition}
\begin{proof}
We refer to  \cite{HJKPZ} for a proof.
\end{proof}
Now we are ready to prove the convergence of spatial configuration.

\begin{theorem}\label{T6.1}
\emph{(Relaxation of spatial configuration)}
Suppose initial data and system parameters satisfy 
\[ \Lambda > \left(\frac{\kappa_3 M_r}{\kappa_2 m_a}\right)^{\beta-\alpha}, \]
where $\Lambda$ is a positive constant defined in \eqref{B-7}, and  let $\{ (W_i,x_i) \}$ be a solution to system \eqref{A-1}. Then, there exists $\{ x_i^\infty \}$ such that
\[
\lim_{t\rightarrow \infty}x_i(t)=x_i^\infty, \quad\forall~ i \in [N].
\]
\end{theorem}
\begin{proof} We consider three cases below. \newline

\noindent $\bullet$~Case A  ($\alpha\neq 2$ and $\beta\neq 2$):~We set the potential function $V$ and the function $f$ by
\begin{align*}
\begin{aligned}
\dot{X}&=-\nabla_XV(X)+f(t),\\
V(X)&:=\frac{\kappa_2}{N}\sum_{i=1}^{N}\sum_{j\neq i} \left((1+\varepsilon_a)\frac{\|x_j-x_i \|^{2-\alpha}}{2-\alpha}\right)-\frac{\kappa_3}{N}\sum_{i=1}^{N}\sum_{j\neq i}\left(\varepsilon_r\frac{\|x_j-x_i \|^{2-\beta}}{2-\beta}\right), \\
f(t)&:=\frac{\kappa_2}{N}\sum_{i=1}^{N}\sum_{j\neq i} \left((\Psi_a^{ij}-1-\varepsilon_a)\frac{x_j-x_i}{\|x_j-x_i\|^\alpha}\right)-\frac{\kappa_3}{N}\sum_{i=1}^{N}\sum_{j\neq i}\left((\Psi_r^{ij} -\varepsilon_r)\frac{x_j-x_i}{\|x_j-x_i\|^\beta}\right).
\end{aligned}
\end{align*}
Note that $X$ is in ${\mathbb{R}^{Nd}}$ when each of $x_i$ is $d-$dimensional. First of all, we claim that the analyticity of the potential function $V$ defined in $\mathbb{R}^{Nd}$. we can restrict the solution space because of the existence of uniform lower bound of distance. We may exclude the hyperplanes $E_{ij}:=\{x\in\mathbb{R}^{Nd}:x_i=x_j\}$ for all $i\neq j$. So we write the restricted space as a union of open connected domain $O_k$
\[
\mathbb{R}^{Nd} \big\backslash\bigcup E_{ij}=\bigcup_k O_k.
\]
By the continuity of the solution, $O_k$ are invariant sets for $X$. Without loss of generality, we may assume $X(0)\in O_1$ and thus $X(t)\in O_1$. That is $O_1$ is an open domain and $V$ is analytic on $O_1$. \newline

Now, we will show that there exists positive constants $C$ and $\lambda$ such that
\[
|f(t)|\leq Ce^{-\lambda t}.
\]
Note that we have verified that there exists a positive constant $c$ such that
\[
\sum_{i=1}^NI_i \leq Ne^{-ct}.
\]
It's clear that for all $i$, 
\[ I_i(t)\leq Ne^{-ct}=:w. \]
In addition, one has 
\begin{align*}
\begin{aligned}
& \Psi_a^{ij} -1-\varepsilon_a=2I_iI_j-I_i-I_j\leq2w^2+w+w\leq(2N^2+2N)e^{-ct}, \\
& \Psi_r^{ij} -\varepsilon_r=I_i+I_j-2I_iI_j\leq2w^2+w+w\leq(2N^2+2N)e^{-ct}.
\end{aligned}
\end{align*}
This means $\Psi_a(W_i,W_j)-1-\varepsilon_a$ and $\Psi_r(W_i,W_j)-\varepsilon_r$ decay in exponential order. Since we have uniform upper and lower bounds of $\|x_i-x_j\|$, there exist positive constants $C$ and $\lambda$ such that
\[
|f(t)|\leq Ce^{-\lambda t}.
\]
Lastly, upper and lower bounds guarantee the boundedness of $V_X(X(t))$ and $\nabla_XV(X(t))$. Similarly, $\frac{d}{dt}|\nabla_XV(X(t))|^2$ is bounded too. Thus we get uniform continuity of $|\nabla_XV(X(t))|^2$ in time. Thus, every hypothesis in Proposition \ref{P4.9} is satisfied to get $X^{\infty}$ that $X(t)$ converges. \newline

\noindent $\bullet$~Case B $(\alpha=2)$:~By setting
\begin{align*}
V(X):=\frac{\kappa_2}{N}\sum_{i=1}^{N}\sum_{j\neq i} \left(2(1+\varepsilon_a)\ln{\|x_j-x_i \|}\right)-\frac{\kappa_3}{N}\sum_{i=1}^{N}\sum_{j\neq i}\left(\varepsilon_r\frac{\|x_j-x_i \|^{2-\beta}}{2-\beta}\right), \end{align*}
we can prove the desired estimate similarly as in Case A. \newline

\noindent $\bullet$ Case C $(\beta=2)$: Similar to Case B, we can show the desired estimate with 
\begin{align*}
V(X)&:=\frac{\kappa_2}{N}\sum_{i=1}^{N}\sum_{j\neq i} \left((1+\varepsilon_a)\frac{\|x_j-x_i \|^{2-\alpha}}{2-\alpha}\right)-\frac{\kappa_3}{N}\sum_{i=1}^{N}\sum_{j\neq i}\left(2\varepsilon_r\ln{\|x_i-x_j \|}\right).
\end{align*}
\end{proof}

\vspace{0.5cm}

Next, we consider a two-particle system with the following initial SIR state:
\begin{align*}
\begin{aligned}
& S_1(0)=s_1,\quad I_1(0)=1-s_1,\quad R_1(0)=0, \\
& S_2(0)=s_2,\quad I_2(0)=1-s_2,\quad R_2(0)=0,
\end{aligned}
\end{align*}
and we set 
\[ s:=\max\{s_1,s_2\}. \]
Consider a system for spatial position:
\begin{align*}
\dot{x}_1=\frac{\kappa_2}{2}\Psi_a^{12} \frac{x_2-x_1}{\|x_1-x_2\|^\alpha}-\frac{\kappa_3}{2}\Psi_r^{12} \frac{x_2-x_1}{\|x_1-x_2\|^\beta},\\
\dot{x}_2=\frac{\kappa_2}{2}\Psi_a^{12} \frac{x_1-x_2}{\|x_1-x_2 \|^\alpha}-\frac{\kappa_3}{2}\Psi_r^{12} \frac{x_1-x_2}{\|x_1-x_2\|^\beta}.
\end{align*}
Note that since $\dot{x}_1+\dot{x}_2=0$, the center of the two particles is constant:
\[  \frac{x_1(t)+x_2(t)}{2} = \frac{x_1^0 +x_2^0}{2}, \quad t \geq 0. \]
 Now, we set 
\[
x :=x_1-x_2,
\]
and derive the equation for $x$:
\[
\dot{x} =\dot{x}_1-\dot{x}_2 =\kappa_2\Psi_a^{12} \frac{-x}{|x|^\alpha}+\kappa_3\Psi_r^{12} \frac{x}{|x|^\beta} =-\frac{x}{|x|^\alpha}\left(\kappa_2\Psi_a^{12} -\kappa_3\Psi_r^{12} \frac{1}{|x|^{\beta-\alpha}}\right).
\]
As long as each coordinate of $x$ is positive,
\[
-\frac{x}{|x|^\alpha}\left((1+\varepsilon_a)-\varepsilon_r\cdot\frac{1}{|x|^{\beta-\alpha}}\right)\leq\dot{x}\leq-\frac{x}{|x|^\alpha}\left(\varepsilon_a-(1+\varepsilon_r)\cdot\frac{1}{|x|^{\beta-\alpha}}\right).
\]
Thus, one has
\[
\min\left\{x^0, \left(\frac{\varepsilon_r}{1+\varepsilon_a}\right)^{\frac{1}{\beta-\alpha}}\right\}<x(t)<\max\left\{x^0, \left(\frac{1+\varepsilon_r}{\varepsilon_a}\right)^{\frac{1}{\beta-\alpha}}\right\},
\]
so the difference $x$ has upper and lower bounds that are uniform in time.

It follows from Lemma $\ref{L5.2}$ that for $\delta_Q=\min\left\{D(0), {\left(\frac{\kappa_3 m_r}{\kappa_2 M_a}\right)}^{\frac{1}{\beta-\alpha}}\right\}$,
\[
\inf_{0\leq t<\infty}D(t)\geq\delta_Q,
\]
and $D(t) = x(t)$. Therefore, $x(t)$ has both upper and lower bounds and by the same reason in Theorem $\ref{T6.1}$, there exists $x^\infty$ such that $x(t)$ converges, as $t$ tends to  infinity. \newline

Now, we consider an equilibrium $(S_i^{\infty}, I_i^{\infty}, R_i^{\infty})$:
\[
\dot{S}_1^\infty=a_{12}S_1^\infty I_2^\infty=0,\quad
\dot{I}_1^\infty=a_{12}S_1^\infty I_2^\infty -bI_1^\infty=0.
\]
Thus we have
\[ I_1^\infty=0, \quad I_2^\infty=0. \] 
This yields
\[ \Psi_a^{12}=1+\varepsilon_a \quad \mbox{and} \quad \Psi_r^{12}=\varepsilon_r. \]
 Since $\dot{x}^\infty=0$,  we have
\[
\kappa_2=\frac{\kappa_3 \varepsilon_r}{(1+\varepsilon_a)|x^\infty|^{\beta-\alpha}}, \quad \mbox{i.e.,} \quad
x^\infty =\left(\frac{\kappa_3\varepsilon_r }{\kappa_2(1+\varepsilon_a)}\right)^{\frac{1}{\beta-\alpha}}.
\]
In addition, we can weaken the condition that guarantees the exponential decay of $I_i$.
\begin{theorem}
Suppose system parameters satisfy
\[
b>\frac{2\kappa_1}{(x^\infty +L)^\gamma}.
\]
and  let $\{(W_i,x_i)\}$ be a solution to system \eqref{B-1}.
Then, there exist positive constants $A, \varepsilon^*$ and $t_0$ such that
\[
\left|I_1(t)+I_2(t)\right| <Ae^{-\frac{b}{2}t},\quad t>t_0,
\]
where $A$ is a positive constant defined by 
\[ A:=e^{t_{\varepsilon^*}}\int_{0}^{t_{\varepsilon^*}}{(I_1(t)+I_2(t)) dt}. \]
\end{theorem}
\begin{proof}
Note that 
\[ \lim_{t\rightarrow\infty}x(t)=x^\infty. \]
Thus, for every $\varepsilon>0$, there exists positive constant $t_\varepsilon$ with respect to $\varepsilon$ such that 
\[  t>t_\varepsilon \quad  \Longrightarrow \quad  |x(t)-x^{\infty}|<\varepsilon. \]
By the given condition, there exists a positive constant $\varepsilon^*$ such that
\[
\frac{\kappa_1}{(x^\infty-\varepsilon^* +L)^\gamma}<\frac{b}{2}.
\]
For $t>t_{\varepsilon^*}$,
\begin{align*}
(I_1(t)+I_2(t))'&=\frac{\kappa_1}{(|x|+L)^\gamma}(S_1I_2+S_2I_1)-b(I_1+I_2) \\
&< \frac{\kappa_1}{(x^\infty-\varepsilon^*+L)^\gamma}(I_2+I_1)-b(I_1+I_2)\\
&< \frac{b}{2}(I_2+I_1)-b(I_1+I_2)= -\frac{b}{2}(I_2+I_1).
\end{align*}
This yields
\[
I_1(t)+I_2(t) < Ae^{-\frac{b}{2}t},\quad t>t_{\varepsilon^*}.
\]
\end{proof}
\begin{remarks} Since
\[
S_i(t)= 1-I_i(t)-R_i(t)=1-I_i(t)-b\int_0^tI_i(\tau)d\tau,
\]
we get
\begin{align*}
\frac{dI_1(t)}{dt}&=a_{12}(t)\left(1-I_1(t)-b\int_0^t{I_1(\tau)d\tau}\right)I_2(t)-bI_1(t),\\
\frac{dI_2(t)}{dt}&=a_{12}(t)\left(1-I_2(t)-b\int_0^t{I_2(\tau)d\tau}\right)I_1(t)-bI_2(t),\\
\frac{dx(t)}{dt}&=-\frac{x}{|x|^\alpha}\left(\kappa_2\Psi_a^{12}-\kappa_3\Psi_r^{12}\cdot\frac{1}{|x|^{\beta-\alpha}}\right).
\end{align*}
The Jacobian matrix of system at $(I_1, I_2) =(0, 0)$ is
\[
\begin{pmatrix}
-b & a_{12}(1-b\int_0^t{I_1(\tau)d\tau}) & 0\\
a_{12}(1-b\int_0^t{I_2(\tau)d\tau}) & -b & 0\\
* & * & *
\end{pmatrix}.
\]
Next, we consider two cases for $b$. \newline

\noindent $\bullet$~Case A:~Consider the case
\begin{align*}
\begin{aligned} 
b &>a_{12}(\infty) \Big[ \Big(1-b\int_0^\infty{I_1(\tau)d\tau} \Big) \Big(1-b\int_0^\infty{I_2(\tau)d\tau} \Big) \Big]^{1/2} \\
   &=\frac{\kappa_1}{(x^{\infty}+L)^\gamma} \Big[ \Big (1-b\int_0^\infty{I_1(\tau)d\tau}  \Big) 
\Big (1-b\int_0^\infty{I_2(\tau)d\tau} \Big) \Big]^{1/2}.
\end{aligned}
\end{align*}
In this case, all eigenvalues are real and eigenvalues for $I_!$ and $I_2$ are negative. This means that $I_1$ and $I_2$ decay to $(0, 0)$ exponentially fast. \newline

\noindent $\bullet$~Case B:~Consider the case
\[ b \leq \frac{\kappa_1}{(x^{\infty}+L)^\gamma} \Big[  \Big (1-b\int_0^\infty{I_1(\tau)d\tau} \Big ) \Big (1-b\int_0^\infty{I_2(\tau)d\tau} \Big) \Big]^{1/2}. \]
The right-hand side is less or equal to 
\[ \frac{\kappa_1}{2(x^{\infty}+L)^\gamma} \Big [ \Big(1-b\int_0^\infty{I_1(\tau)d\tau} \Big)+ \Big(1-b\int_0^\infty{I_2(\tau)d\tau} \Big) \Big]. \]
This yields
\[
\int_0^\infty{(I_1+I_2)(t)dt} \leq \frac{2}{b}-\frac{2(x^\infty +L)^\gamma}{\kappa_1}.
\]
Note that the left-hand side $\int_0^\infty{(I_1+I_2)(t)dt}$ denotes the number of  total infections. Thus, we can say that $I_i(t)$ decays to zero exponentially fast, or the total number of 
 total infection is bounded. 

\end{remarks}

\section{Numerical Simulations}\label{sec:7}
In this section, we provide several numerical examples to confirm analytical convergence results that we have shown in previous sections. At first, we set initial condition by $(S_i(0),I_i(0),R_i(0))=(1,0,0)$ for uninfected 16 people and $(S_i(0),I_i(0),R_i(0))=(0.1,0.9,0)$ for infected 4 people. Initial location is set by random seed in 3 by 3 plane. We set 
\[
\kappa_1=1,\quad \varepsilon_a=\varepsilon_r=0.2,
\]
and we use the fourth order Runge-Kutta method for all simulations.  \newline

\noindent $\bullet$~(Convergence of $(S_i, I_i, R_i)$): In Figure 1 and Figure 2, we used 
\[
b=0.4,\quad \gamma=1,\quad L=1,
\]
here. In Theorem \ref{T4.1},  we have shown that $S_i$, $I_i$, and $R_i$ converges. To authenticate this assertion by numerical method, we observed convergence of $S_i$, $I_i$, and $R_i$.
\begin{figure}[h]
\centering
\subfigure[$S,I,R$ for $\alpha=1$, $\beta=2$, $\kappa_2=1$, $\kappa_3=5$]{\includegraphics[height=5cm]{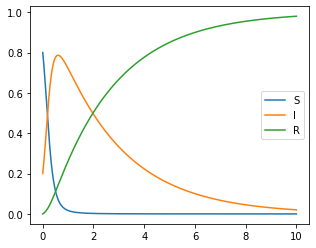}}
\subfigure[$S,I,R$ for $\alpha=2$, $\beta=5$, $\kappa_2=10$, $\kappa_3=1$]{\includegraphics[height=5cm]{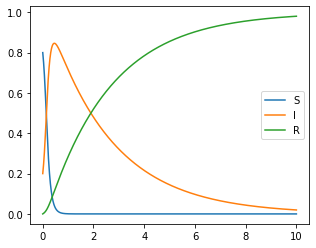}}
\caption{Convergence of $S,I$ and $R$}
\end{figure}
Moreover, by observing $\dot{S}_i$, $\dot{I}_i$, and $\dot{R}_i$ that converges to zero, we obtain that $S_i,I_i$ and $R_i$ converge.
\begin{figure}[h]
\centering
\subfigure[$\alpha=1$, $\beta=2$, $\kappa_2=1$, $\kappa_3=5$]{\includegraphics[height=5cm]{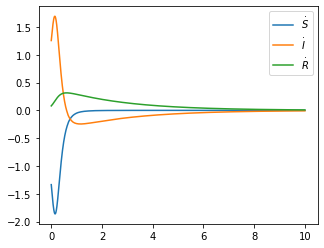}}
\subfigure[$\alpha=2$, $\beta=5$, $\kappa_2=10$, $\kappa_3=1$]{\includegraphics[height=5cm]{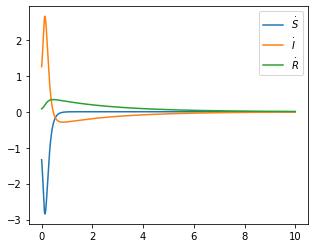}}
\caption{Convergence of $\dot{S},\dot{I}$ and $\dot{R}$}
\end{figure}

\vspace{0.2cm}

\noindent $\bullet$~(Exponential decay of $\sum_i{I_i}$):~For arbitrary $N$, we showed that 
$\sum_i{I_i}$ decays to zero exponentially fast when $b>\frac{\kappa_1(N-1)}{L^\gamma}$ as in Theorem \ref{T4.1}. In Figure 3 and Figure 4, we set 
\[ b=1,~ \gamma=3,~ L=3 \quad \mbox{so that} \quad \frac{\kappa_1(N-1)}{L^\gamma}\approx 0.7037 < 1 = b, \]
and observed two simulations based on two set of system parameters:
\[ (\alpha, \beta, \kappa_2, \kappa_3) = (1,2,1,5), \quad (2,5, 10, 1). \]
\begin{figure}[h]
\centering
\subfigure[$S,I,R$ for $\alpha=1$, $\beta=2$, $\kappa_2=1$, $\kappa_3=5$]{\includegraphics[height=5cm]{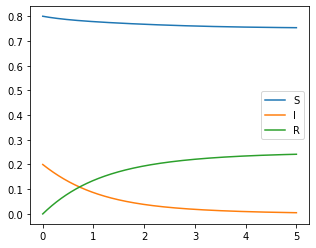}}
\subfigure[$\log\sum_i{I_i}$ for $\alpha=1$, $\beta=2$, $\kappa_2=1$, $\kappa_3=5$]{\includegraphics[height=5cm]{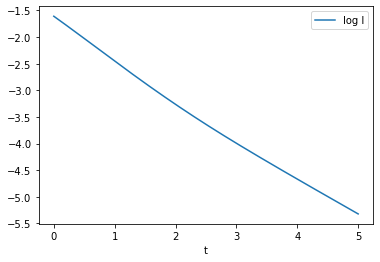}}
\caption{Exponential decay of $\log\sum_i{I_i}$}
\end{figure}
\begin{figure}[h]
\centering
\subfigure[$S,I,R$ for $\alpha=2$, $\beta=5$, $\kappa_2=10$, $\kappa_3=1$]{\includegraphics[height=5cm]{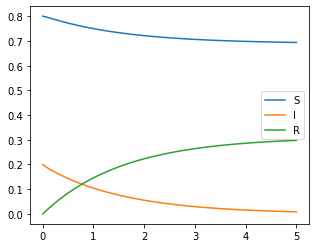}}
\subfigure[$\log\sum_i{I_i}$ for $\alpha=2$, $\beta=5$, $\kappa_2=10$, $\kappa_3=1$]{\includegraphics[height=5cm]{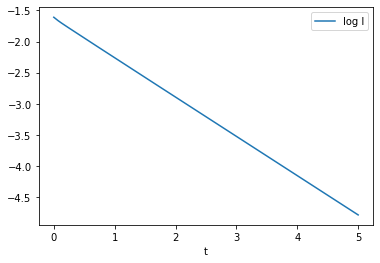}}
\caption{Exponential decay of $\log \Big( \sum_i{I_i} \Big)$}
\end{figure}

%

\vspace{0.1cm}

\noindent $\bullet$~(Effect of social distancing): Since $\kappa_2$ and $\kappa_3$ are coefficients of attracting and repulsion forces, respectively,  larger ratio $\frac{\kappa_3}{\kappa_2}$ represents for more intensive social distancing. We observe that the maximal value of $\sum_i{I_i}$ decreases, as $\frac{\kappa_3}{\kappa_2}$ increases. We set $L=1$ and $\kappa_2=1$, and changed values of $b$ and $\Psi$ with $\kappa_3=1,5,10,50$. We plot graph for the cases of $\frac{\kappa_3}{\kappa_2}=1,5,10$ and $50$.
\begin{figure}[h]
\centering
\subfigure[$\frac{1}{N}\sum_i{I_i}$ for $\beta=0.2$, $\Psi=1$]{\includegraphics[height=5cm]{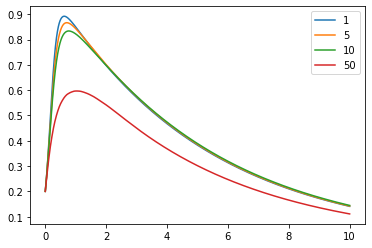}}
\subfigure[$\frac{1}{N}\sum_i{I_i}$ for $\beta=0.4$, $\Psi=3$]{\includegraphics[height=5cm]{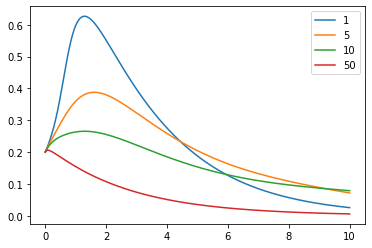}}
\caption{Average of ${I_i}$, as $\frac{\kappa_3}{\kappa_2}$ changes}
\end{figure}

\vspace{0.2cm}

\begin{figure}[h]
\includegraphics[height=10cm]{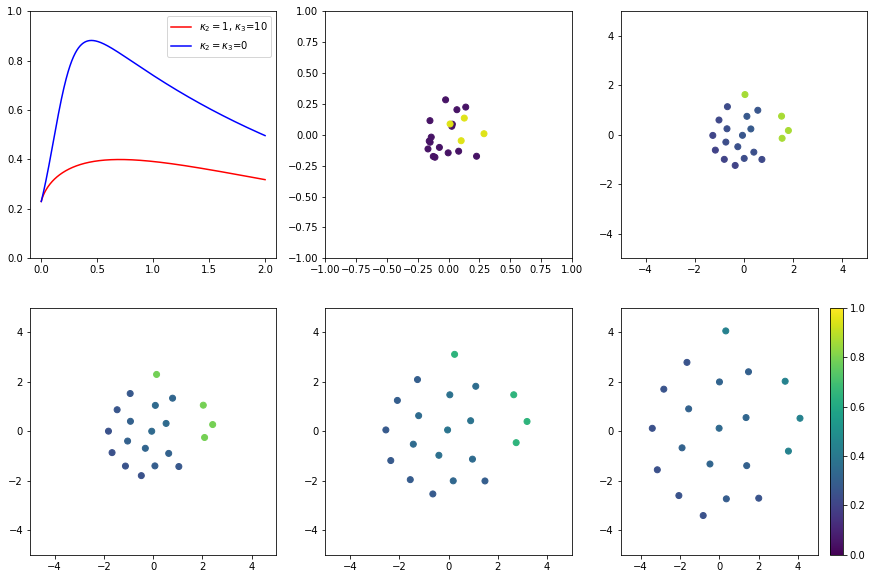}
\caption{Behavior of particles in $b=0.4, \Psi=3, L=1$}
\label{move1}
\end{figure}

\noindent $\bullet$~(Behavior of the particles): We observed that infectious particles aggregated and they were isolated from non-infectious factors. It implies that this model can effectively reduce the number of infectious particles by adjusting the coefficients. We illustrated the behavior of the particles in
\[
b=0.4,~ \gamma=3,~ L=1,
\]
and $\kappa_2=1$, $\kappa_3=10$. The purple ones represents for non-infectious particles. Moreover, we compared with the model that has no attracting and no repulsing term by setting $\kappa_2=\kappa_3=0$. It is illustrated in the first graph in Figure \ref{move1} , red one is for $\kappa_2=1$, $\kappa_3=10$ and blue one is for  $\kappa_2=\kappa_3=0$. We could decrease the average of $I_i(t)$.
We set
\[
b=0.2,~ \gamma=3,~ L=3,~ \alpha=2,~ \beta=5
\]
and repeated the process.
\begin{figure}[h]
\includegraphics[height=10cm]{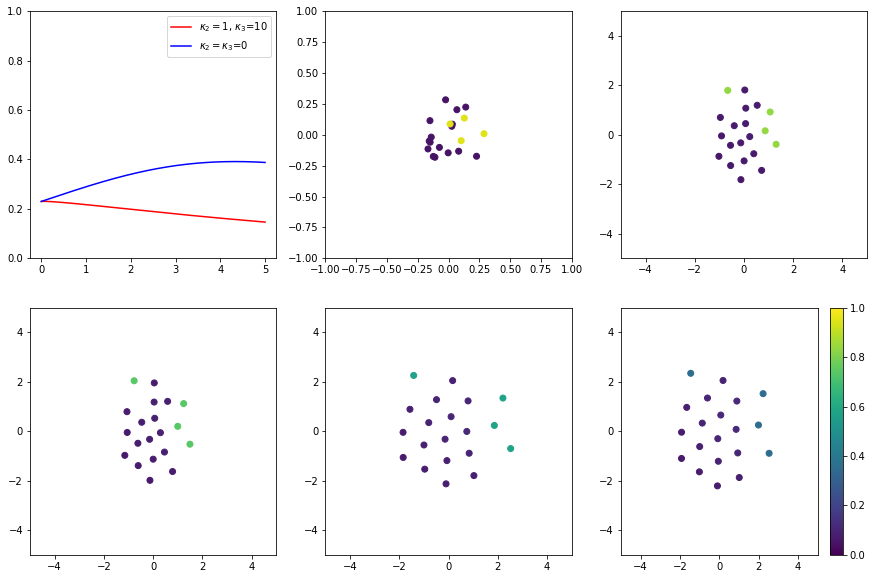}
\caption{Behavior of particles in $b=0.2, \Psi=3, L=3$}
\label{move2}
\end{figure}

\section{conclusion} \label{sec:8}
In this paper, we have studied the emergent behaviors of a flock with SIR internal states, and have presented a new particle model with an aggregate property and epidemic internal forces by combining the SIR model and the swarmalator model. We considered that each particle has a state vector which is a kind of internal state following the SIR dynamics. From this argument, we could obtain the SIR-aggregation model by imposing the repulsive/attractive force between two particles depending on the distance between two particles and the internal states. We imposed strong repulsive force if one of them has a high probability of infection to model the social distancing, in contrast, we imposed weak repulsive force if both of them has low probabilities of infection. From this modeling, we modeled the social distancing of particles. We expect that we can do numeric experiments to find the efficiency of social distance for given strategies. We also provided the theoretical result of the SIR model, for example, nonexistence of finite-time collisions, positive lower bound for minimal relative distance, and a uniform upper bound for spatial diameter. From the numerical simulations, we could check that the isolation of the infected particle is an effective way to reduce the number of infected particles. Since we only considered that the recovered particles were never infected again, we have monotonicity on the number of susceptible/recovered particles. From those properties, we could prove the emergent dynamics. In our proposed model, we did not consider the reinfection which can happen in reality. Thus, we will leave this issue for a future work.

\end{document}